\documentclass[10pt,reqno]{amsart}

\usepackage{amssymb, amsthm, amsmath, latexsym}
\usepackage{amsfonts}
\usepackage[dvips]{graphicx}
\usepackage{psfrag}

\usepackage{color,mathrsfs,epsfig}
\newcommand{\bp}{\begin{proof}}
\newcommand{\ep}{\end{proof}}
\newcommand{\be}{\begin{equation}}
\newcommand{\ee}{\end{equation}}

\newcommand{\e}{\varepsilon}

\newcommand{\ue}{U^{\varepsilon}}

\newcommand{\R}{\mathbb R}
\newcommand{\N}{\mathbb N}

\newcommand{\unpo}{\mathcal{O}(1)}

\newcommand{\tchar}{\mathfrak T_k}
\newtheorem{thm}{Theorem}[section]
\newtheorem{lemma}[thm]{Lemma}

\newtheorem{hyp}{Hypothesis}

\numberwithin{equation}{section}

\theoremstyle{definition}
\newtheorem{remark}{Remark}[section]
\renewenvironment{proof}[1][Proof]{\textbf{#1.} }{\ \rule{0.5em}{0.5em}}

\numberwithin{equation}{section}
\begin{document}

\bibliographystyle{plain}
\title[Uniqueness criterion for the boundary Riemann problem]{A uniqueness criterion for viscous limits of boundary Riemann problems}
\author{Cleopatra Christoforou \and Laura V.  Spinolo}

\address{C.C.: Department of Mathematics\\
University of Houston\\
4800 Calhoun Road\\
Houston, TX 77204, USA 
and 
Department of Mathematics and Statistics\\
University of Cyprus\\
1678 Nicosia, Cyprus}
\email{Christoforou.Cleopatra@ucy.ac.cy}
\urladdr{http://192.42.1.1/$\sim$kleopatr}

\address{L.V.S.: Centro De Giorgi \\
Collegio Puteano, Scuola Normale Superiore \\
Piazza dei Cavalieri 3 \\
56126 Pisa, Italy}
\email{laura.spinolo@sns.it }
\urladdr{http://www.math.northwestern.edu/$\sim$spinolo}
\maketitle
\begin{abstract}
We deal with initial-boundary value problems for systems of conservation laws in one space dimension and we focus on the boundary Riemann problem. It is known that, in general, different viscous approximations provide different limits. In this paper, we establish sufficient conditions to conclude that two different approximations lead to the same limit. 
As an application of this result, we show that, under reasonable assumptions, the self-similar second-order approximation 
$$ 
    \partial_t U^{\varepsilon  }+  \partial_x   F( U^{\varepsilon})=\e\, t \,  \partial_x \big( 
               B(U^{\varepsilon}) \partial_x U^{\varepsilon}  \big)
$$
and the classical viscous approximation 
$$
\partial_t U^{\varepsilon  }+  \partial_x   F( U^{\varepsilon})=\e  \partial_x \big( 
               B(U^{\varepsilon}) \partial_x U^{\varepsilon}  \big) 
$$               
provide the same limit as $\varepsilon \to 0^+$. 
Our analysis applies to both the characteristic  and the non characteristic case. We require neither genuine nonlinearity nor linear degeneracy of the characteristic fields.

Date: \today 
\end{abstract}
\section{Introduction and main results}
\label{s:intro}
We are interested in the system of conservation laws 
\be\label{S1: hyp}
	\partial_t U+\partial_x F(U)=0, 
\ee
where the unknown function $U(t, x)$ takes values in $\R^n$, the variables $t$ and $x$ are both scalar and the flux $F: \R^n \to \R^n$ is a smooth function. We assume that system~\eqref{S1: hyp} is strictly hyperbolic, namely for every $U \in \R^n$ the Jacobian matrix $DF(U)$ has $n$ real and distinct eigenvalues
$$
    \lambda_1 (U) < \lambda_2 (U) < \dots < \lambda_n (U).
$$
For an exposition of the current state of the theory of systems of conservation laws, we refer to the books by Dafermos~\cite{D} and Serre~\cite{Serre:book}. 

In this paper, we study initial boundary value problems and we establish conditions ensuring that two different viscous approximations lead to the same solution in the limit. Then, we apply our result to self-similar second-order approximations and the classical viscous approximations.

To contextualize the problem, let us first consider the Cauchy problem obtained by coupling~\eqref{S1: hyp} with the initial datum
\be
\label{e:cauchy}
       U(0, x) = U_c(x). 
\ee
The so-called \emph{Riemann problem} is posed in the case that $U_c(x)$ has the form 
\be 
 \label{e:i:rie}
       U_c( x) = 
       \left\{
       \begin{array}{ll}
                  U^+ & \text{if $x>0$ } \\
                  U^- & \text{if $x< 0$,  } \\  
       \end{array}
       \right. 
\ee      
where $U^-, \; U^+ $ are two given constant states in $\R^n$. The Riemann problem~\eqref{S1: hyp},~\eqref{e:i:rie} has been extensively studied and its analysis provides important information on both the local (in space-time) and the long-time behavior of the solution of a general Cauchy problem~\eqref{S1: hyp}-\eqref{e:cauchy}. Also, it serves as a building block for the construction of different approximation schemes (e.g. the random choice method of Glimm~\cite{Glimm}, the wave front-tracking algorithm, see Bressan et al~\cite{Bre:book} and Holden and Risebro~\cite{HoldenRisebro}) that lead to existence and uniqueness results for general Cauchy problems. It should be noted that these results hold for general flux $F$ under the hypothesis that the total variation of the Cauchy datum $U_c(x)$ is small enough. Suitable counterexamples ensure that handling data with large or unbounded total variation requires the assumption of more restrictive conditions on the structure of the flux $F$ in~\eqref{S1: hyp}. In the present paper, we deal with general flux functions and, therefore, we focus on Riemann data $U^+$, $U^-$ in~\eqref{e:i:rie} that are sufficiently close.   

One of the major difficulties posed by the Riemann problem is the non-uniqueness of distributional solutions. Various selection principles, often motivated by physical considerations, have been introduced in an attempt at singling out a unique solution: see, for example, the entropy admissibility criterion and the conditions named after Lax~\cite{lax} and Liu~\cite{Liu:TAMS,Liu:adm}.  A solution of the Riemann problem~\eqref{S1: hyp},~\eqref{e:i:rie} satisfying suitable admissibility criteria was first constructed by Lax~\cite{lax},  Liu~\cite{Liu:rie} and later by Tzavaras~\cite{Tz} and Bianchini~\cite{Bia:riemann}. 

In particular, in~\cite{Bia:riemann}, solutions of~\eqref{S1: hyp},~\eqref{e:i:rie} are obtained by considering the viscous approximation
\be
\label{S1: viscous}         
               \partial_t U^{\varepsilon  }+  \partial_x   F( U^{\varepsilon})=\e\; \partial_x \big( 
               B(U^{\varepsilon}) \partial_x U^{\varepsilon}  \big) \;,
\ee
where $B$ is an $n \times n$ matrix which depends on the physical model under consideration and $\varepsilon$ is a positive parameter. As $\varepsilon \to 0^+$, the family of functions $U^{\varepsilon}$ is expected to converge to a distributional solution of~\eqref{S1: hyp}. It is a challenging open problem to establish a rigorous proof of the convergence $\ue\to U$ in the general case, but results have been achieved in specific cases. See Bianchini and Bressan~\cite{BiaBrevv} and the references in the books~\cite{D, Serre:book}. It should also be noted that the analysis in~\cite{Bia:riemann} ensures that the function constructed in~\cite{Bia:riemann, lax,Liu:rie,Tz} is the unique solution of~\eqref{S1: hyp},~\eqref{e:i:rie} which is self-similar and can be obtained by patching together a countable number of rarefaction waves and of shocks (or contact discontinuities) satisfying the admissibility condition named after Liu. As a consequence, the limit as $\e\to0^+$ of the family of functions $\ue$ in~\eqref{S1: viscous} does not depend on the choice of $B$. As we see in the following, this is not true in the case of initial-boundary value problems.

Another approach to the analysis of the Riemann problem was introduced independently by Dafermos~\cite{D1}, Kalasnikov~\cite{kal} and Tupciev~\cite{Tu1966} and is based on the analysis of the family of parabolic problems   
\be
\label{S1: t-viscous}         
               \partial_t U^{\varepsilon  }+  \partial_x   F( U^{\varepsilon})=\e\, t \,  \partial_x \big( 
               B(U^{\varepsilon}) \partial_x U^{\varepsilon}  \big),
\ee
coupled with the initial datum
\be 
 \label{e:vappd}
      \ue(0, x) = 
       \left\{
       \begin{array}{ll}
                  U^+ & \text{if $x>0$ } \\
                  U^- & \text{if $x< 0$.  } \\  
       \end{array}
       \right. 
\ee 
Because of the presence of the ``$t$'' factor in the second order term, the Cauchy problem~\eqref{S1: t-viscous}-\eqref{e:vappd} admits self-similar solutions in the form $U^{\varepsilon  } (t, x) = V^{\varepsilon  } (x /t)$. Compactness results for the solutions of~\eqref{S1: t-viscous}-\eqref{e:vappd} have been established under suitable hypotheses on $B$. See in particular Tzavaras~\cite{Tzavaras:JDE,Tz}, Joseph and LeFloch~\cite{JoLeF:nc} and the references in Dafermos~\cite[Section 9.8]{D}.

In the present paper, we consider the initial-boundary value problem for the system of conservation laws~\eqref{S1: hyp}, hence we assume that $x\geq 0$. The analogous to the Riemann problem is the so-called \emph{boundary Riemann problem}, obtained by coupling~\eqref{S1: hyp} with the Cauchy and Dirichlet data
\be
\label{e:cdrie}
               U(0, x) = U_0, \quad \text{for}\; x>0,\qquad U(t, 0) = \bar U,\quad \text{for}\; t>0, 
\ee
where $U_0$ and $\bar U$ are two given constant states. Since we are interested in solutions of small total variation, we restrict to the case that $|\bar U - U_0|$ is small. The analysis of the boundary Riemann problem provides again information on the local and the long-time behavior of the solution of a general initial-boundary value problem and, moreover, it constitutes the building block for the construction of approximation schemes that lead to general existence and uniqueness results. In addition, because of the presence of the boundary, additional challenges are raised compared to the Riemann problem. First of all, the initial-boundary value problem~\eqref{S1: hyp},~\eqref{e:cdrie} is in general ill-posed, i.e. it posseses no solution. Therefore, a notion of admissible set of boundary values can be introduced, see Dubois and LeFloch~\cite{DubLeF}.    

Other challenges arise when studying the viscous approximation 
\be
\label{e:brs}
\left\{   
\begin{array}{ll}      
               \partial_t U^{\varepsilon  }+  \partial_x   F( U^{\varepsilon})=\e\; \partial_x \big( 
               B(U^{\varepsilon}) \partial_x U^{\varepsilon}  \big)  \\
                 U^{\varepsilon}(0, x) = U_0, \qquad x>0\\
                 U^{\varepsilon}(t, 0) = U_D, \qquad t>0.
\end{array}
\right.
\ee
It should be noted that, if the matrix $B$ is invertible, the above initial-boundary value problem is well-posed for every $U_D$, while if $B$ is singular, one has to introduce a more complicated formulation of the boundary condition, see for example Bianchini and Spinolo~\cite{BiaSpi}. We emphasize that $B$ is singular in most of the physically relevant examples, but to simplify the exposition, we focus on the case of $B$ invertible. 

The first difficulty one has to address when studying~\eqref{e:brs} is the following: assume that as $\varepsilon \to 0^+$ the family $U^{\varepsilon}$ converges, in a suitable topology, to some limit $U$. Also, assume that $U$ is self-similar function, $U(t, x) = V(x/t)$ and that $V$ has bounded total variation, so that $\lim_{\xi \to 0^+} V(\xi)$ is well defined. Because of boundary layer phenomena, in general
$$
    U(t,0) = \lim_{\xi \to 0^+} V(\xi)  \neq U_D. 
$$     
Also, as pointed out by Gisclon and Serre~\cite{Gis,GisSerre}, in general, if $F$, $U_0$ and $U_D$ in~\eqref{e:brs} are fixed and the matrix $B$ varies, then the limit $U$ varies. To understand this behavior, let us first focus on the case that the boundary is \emph{non characteristic}, namely all the eigenvalues of the Jacobian matrix $DF(U)$ are bounded away from $0$ and let us denote the trace by $$\bar U\doteq\lim_{\xi \to 0^+} V(\xi).$$ See condition~\eqref{e:nchar} in Section~\ref{sss:hyp} for the rigorous definition of non characteristic boundary. Then one expects that there is a boundary layer of~\eqref{S1: viscous} connecting $\bar U$ with $U_D$, namely the system 
\be 
\label{e:i:blnc}
\left\{   
\begin{array}{ll}      
              B(W) W' = F(W) - F(\bar U)  \\
              W(0) = U_D  \qquad \lim_{y \to + \infty} W(y) = \bar U \\
\end{array}
\right.
\ee
admits a solution $W(y)$. Here, $'$ denotes the first derivative. Since the solvability of~\eqref{e:i:blnc} depends on $B$, then also $\bar U$ does. The boundary characteristic case occurs if an eigenvalue of the Jacobian matrix can attain the value $0$, see conditions~\eqref{e:c}-\eqref{e:lk} in Section~\ref{sss:hyp} for the precise definition. The boundary case is more complicated to handle because of the possible interactions between boundary layers, traveling waves and contact discontinuities, see the analysis in Ancona and Bianchini~\cite{AnBia} and Bianchini and Spinolo~\cite{BiaSpi}. However, one again obtains that in general the limit depends on the viscosity matrix $B$. 

In~\cite{JosephLeFloch}, Joseph and LeFloch studied the viscous approximation 
\be
\label{e:ssbrs}
\left\{   
\begin{array}{ll}      
               \partial_t U^{\varepsilon  }+  \partial_x   F( U^{\varepsilon})=\e \, t \,  \partial_x \big( 
               B(U^{\varepsilon}) \partial_x U^{\varepsilon}  \big)  \\
                 U^{\varepsilon}(0, x) = U_0, \qquad x>0,\\
                 U^{\varepsilon}(t, 0) = U_D,\qquad t>0,
\end{array}
\right.
\ee
in the case that $B$ is the identity matrix, $B(U) \equiv I$. They established compactness results and provided a detailed description of the limit in both the case of a non characteristic and (under some technical assumptions) a characteristic boundary. In~\cite{JoLeF:nc}, among other results, Joseph and LeFloch extended the analysis to more general viscosity matrices. 

The motivation of this work is to investigate whether the limit of the  classical viscous approximation~\eqref{e:brs} and of the self-similar viscous approximation~\eqref{e:ssbrs} coincide. In the case of the Riemann problem, the limit of~\eqref{S1: viscous},~\eqref{e:vappd} and of~\eqref{S1: t-viscous}-\eqref{e:vappd} are expected to coincide, because of the uniqueness result in Bianchini~\cite{Bia:riemann}. However, the case of the boundary Riemann problem is a priori much less clear, because of the work by Gisclon and Serre~\cite{Gis,GisSerre} stating that, in general, the limit of the viscous approximation depends on the viscosity matrix. On the other hand, the analysis in Joseph and LeFloch~\cite{JosephLeFloch, JoLeF:nc} ensures that the equation of the boundary layers of the self-similar approximation~\eqref{S1: t-viscous} is the same as the equation of the boundary layers of the classical viscous approximation~\eqref{S1: viscous}, namely it is  
\be 
\label{e:bounlay}
       \big[ B(U) U' \big]' = \big[ F(U) \big]'\;.
\ee
Hence, to determine whether the limits of~\eqref{e:brs} and~\eqref{e:ssbrs} coincide, a first step is to investigate whether there exists a unique solution of the limiting system which can be constructed by patching together boundary layers satisfying~\eqref{e:bounlay} and a countable number of rarefaction waves and shocks (or contact discontinuities) that are admissible in the sense of Liu. 

Our first result provides an affirmative answer to the previous question. Here, we state the main result of this paper:
 
\begin{thm}
\label{t:ui} Assume that Hypotheses $1$--$3$ given in Section~\ref{sss:hyp} hold, let $U_0, U_D$ be two given constant states in $ \R^n$ and consider the system of conservation laws 
\be
\label{e:t:cl}
    \partial_t U + \partial_x F(U)=0.
\ee
Then, there exist positive constants $\delta$ and $C$, $\delta$ sufficiently small, such that, if $|U_0 - U_D| \leq \delta$, then there exists a unique distributional solution $U$ of~\eqref{e:t:cl} satisfying the following properties:
\begin{enumerate}
\item the function $U$ attains the Cauchy datum $U(0, x) = U_0$ for a.e. $x$;
\item the function $U$ is self-similar, namely $U(t, x) = V(x/ t)$ for a suitable function $V$ which also satisfies $\mathrm{TotVar} \, V \leq C \delta$;
\item all the shocks and the contact discontinuities of $U$ are admissible in the sense of Liu;
\item let $I$ be any open interval where $V$ is continuous, then $V$ is differentiable in $I$;  
\item there exists a value $\underline U$ with $|U_0 - \underline U| \leq C \delta$, such that the following two conditions hold.
\begin{enumerate}
\item If we denote by $\bar U$ the trace of $U$ on the t-axis, namely
\be
\label{e:htrace}
       \bar U : = \lim_{ \xi \to 0^+} V(\xi ),
\ee
then $F(\underline U) = F(\bar U)$ and the shock or the contact discontinuity connecting 
$\underline U$ (on the left) with $\bar U$ (on the right) is admissible in the sense of Liu.  
\item There exists a boundary layer $W$ satisfying 
\begin{equation}
\label{e:i:bl}
\left\{
\begin{array}{ll}
          B(W) W' = F(W) - F(\underline U) \\
          W(0) = U_D \qquad \lim_{y \to + \infty } W (y) = \underline U
\end{array}
\right.
\ee
and $|W(y)-\underline U|\le C\delta$, $|W'(y)| \leq C \delta$ for every $y$.
\end{enumerate}
\end{enumerate}
\end{thm}

It should be mentioned that Liu's admissibility condition is stated in Section~\ref{sus:rie}, see formula~\eqref{e:liu}. For an extended discussion, see also Dafermos~\cite[Chapter 8.4]{D}. 

Some remarks are here in order. First, the novelty in Theorem~\ref{t:ui} is the uniqueness part, since the existence part is already known. Second, if the boundary is non-characteristic then $F(\underline U) = F (\bar U)$ implies that $\underline U = \bar U$, hence~\eqref{e:i:bl} reduces to~\eqref{e:i:blnc} and the proof of Theorem~\ref{t:ui} is easier.
 In the boundary characteristic case, one has to handle the possibility of a zero-speed shock (or contact discontinuity), so the value $\underline U$ comes into play and the analysis is considerably more complicated. Last, the hypotheses introduced in Section~\ref{sss:hyp} are quite standard in this context and they imply that the viscosity matrix $B$ is strictly stable in the sense of Majda and Pego~\cite{MajdaPego} (see Corollary $2.2$ therein) and they are the same assumptions considered by Gisclon and Serre in~\cite{Gis,GisSerre}. In particular, we require neither genuine nonlinearity nor linear degeneracy of the characteristic fields. Our hypotheses imply that the viscosity matrix $B$ is invertible. The extension to the case of a singular viscosity matrix, which is more interesting in view of physical applications, raises no apparent difficulty provided that the assumptions introduced by Kawashima and Shizuta~\cite{KawShi:normal} and a condition of so-called \emph{block linear degeneracy} introduced in Bianchini and Spinolo~\cite{BiaSpi} are all satisfied. However, in the present paper, to simplify the exposition, we restrict to the case of an invertible viscosity matrix.  

An application of Theorem~\ref{t:ui} to the classical and the self-similar viscous approximation shows that the two limits coincide even in the case of initial-boundary value problems. Here is the result:

\begin{thm}
\label{t:vl}
Let Hypotheses $1$--$3$ in Section~\ref{sss:hyp} hold, then there exist constants $C$ and $\delta$, $\delta$ sufficiently small, such that, given two values $U_0$, $U_D \in \R^n$ satisfying $|U_0 - U_D| \leq \delta$, the following holds. Consider the classical viscous approximation 
\be
\label{e:cvl}
    \left\{
    \begin{array}{ll}
               \partial_t U^{\varepsilon  }+  \partial_x   F( U^{\varepsilon})=\e\; \partial_x \big( 
               B(U^{\varepsilon}) \partial_x U^{\varepsilon}  \big) \\
               U^{\varepsilon}(0, x) = U_0, \qquad x>0,\\
                U^{\varepsilon} (t, 0) = U_D, \qquad t>0,
         \end{array}
    \right.
\ee
and the ``self-similar" approximation 
\be
\label{e:ssvl}
    \left\{
    \begin{array}{ll}
               \partial_t Z^{\varepsilon  }+  \partial_x   F( Z^{\varepsilon})=\e\; t\;\partial_x \big( 
               B(Z^{\varepsilon}) \partial_x Z^{\varepsilon}  \big) \\
               Z^{\varepsilon}(0, x) = U_0 ,\qquad x>0,\\
               Z^{\varepsilon} (t, 0) = U_D, \qquad t>0,
         \end{array}
    \right.
\ee
and assume that  as $\e\to0^+$
$$
    U^{\varepsilon} \to U \;\; \text{and} \; \;
    Z^{\varepsilon} \to Z \quad \text{in $L^1_{\mathrm{loc}} \big( [0, + \infty[ \times [0, + \infty[  \big)$}.
$$
If both $U$ and $Z$ satisfy conditions (1)-(5) in the statement of Theorem~\ref{t:ui}, then $U(t, x) = Z(t, x)$ for a.e.~$(t, x)\in[0,+\infty[\times[0,+\infty[$. 
\end{thm}
 As a final remark, we point out that assuming that both the limits $U$ and $Z$ satisfy conditions (1)-(5) in the statement of Theorem~\ref{t:ui} is reasonable. Indeed, concerning the classical viscous approximation~\eqref{e:cvl}, these hypotheses are satisfied in the case when the global in time convergence has been proved, namely when $B (U) \equiv I$, see  Ancona and Bianchini~\cite{AnBia} and also the analysis in Bianchini and Bressan~\cite[Section 14]{BiaBrevv}. As pointed out before, concerning the self-similar viscous approximation~\eqref{e:ssvl}, convergence results have been established by Joseph and LeFloch~\cite{JosephLeFloch, JoLeF:nc} under suitable assumptions on the matrix $B$: their analysis shows that, under the same assumptions, condition ($5$b) holds. Conditions (1)-(4) in the statement of Theorem~\ref{t:ui} are reasonable in view of the analysis in Tzavaras~\cite[Sections 8-9]{Tz} and Dafermos~\cite[Section 9.8]{D}. 

The paper is organized as follows: the hypotheses exploited in the work are discussed in Section~\ref{sss:hyp}, whereas in Section \ref{s:pre} we present some preliminary results that are used in the following section. The proof of Theorem~\ref{t:ui} is established in Section \ref{s:uni} and Theorem~\ref{t:vl} is an immediate consequence of Theorem~\ref{t:ui}.

\subsection{Hypotheses}
\label{sss:hyp}
Here, we introduce the hypotheses of Theorem~\ref{t:ui}. The first one is the standard assumption of strict hyperbolicity.
\begin{hyp}
\label{h:hyp}
          The flux $F: \R^n \to \R^n$ is smooth and, for every $U \in \R^n$, the Jacobian matrix $DF (U)$ 
          has $n$ real and distinct eigenvalues, 
         \begin{equation}
         \label{e:eig}
             \lambda_1 (U) < \lambda_2 (U) < \dots < \lambda_n (U).
         \end{equation}
 \end{hyp} 
Next, we also assume that the system of conservation laws is endowed with a strictly convex entropy.
\begin{hyp}
\label{h:entr} There are functions $\eta: \R^n \to \R$ and $q: \R^n \to \R$ such that
         \be 
         \label{e:entro}
             \nabla \eta (U) \cdot DF (U) = \nabla q (U) 
         \ee
         and $D^2 \eta (U)$ is strictly positive definite for every $U \in \R^n$. Here,  $``\cdot"$ denotes the standard dot product.
\end{hyp}

Last, we assume that the viscosity matrix $B$ appearing in~\eqref{S1: viscous} and~\eqref{S1: t-viscous} is dissipative.
\begin{hyp}
\label{h:diss}
        For any given compact set $H \subseteq \R^n$, there exists a constant $\alpha_H >0$ such that, for every $U \in H$, 
        \be 
        \label{e:diss}
            D^2 \eta (U)  \vec v \cdot B(U) \vec v \ge \alpha_H |\vec v|^2 \qquad\forall \, \vec v \in \R^n.
        \ee
\end{hyp}
We also quote Proposition $3$ in Gisclon~\cite{Gis}, since it is needed in Section~\ref{s:uni}.
\begin{lemma}
\label{l:same}
       If Hypotheses~\ref{h:hyp}--\ref{h:diss} hold, then, for every $U \in \R^n$, the number of eigenvalues of $DF(U)$ with positive (resp. negative) real part is the same as the number of eigenvalues of $B^{-1}(U) DF(U)$ with positive (resp. negative) real part.
\end{lemma}
Now, we introduce the following standard terminology. The boundary in problem~\eqref{S1: hyp},~\eqref{e:cdrie} is non characteristic if all the eigenvalues of the Jacobian matrix $DF(U)$ are bounded away from zero, namely if there exists a constant $c>0$ such that 
\be
\label{e:nchar}
    \lambda_1 (U) < \dots < \lambda_{n-p} (U)< -c <0 < c < \lambda_{n-p+1}(U) < \dots < \lambda_n (U)   
\ee
for every $U \in \R^n$. It is actually sufficient to assume that~\eqref{e:nchar} holds for every $U$ such that $|U-U_0|\le 3\,C\delta$, where $C$ and $\delta$ are the same constants as in the statement of Theorems~\ref{t:ui} and~\ref{t:vl}.

Conversely, the boundary is characteristic if~\eqref{e:nchar} is violated, namely if one eigenvalue of $DF(U)$ can attain the value $0$. If the boundary is characteristic, there are positive constants $K, c>0$ satisfying
\be
\label{e:c}
    \lambda_1 (U) < \dots < \lambda_{k-1}(U) < -c <0 < c < \lambda_{k+1} (U)< \dots < \lambda_n (U)   
\ee
and 
\begin{equation}
\label{e:lk}
   | \lambda_k (U) | \leq K \delta
\end{equation}
for every value of $U$ such that $|U-U_0|\le 3C\delta$, where $C$ and $\delta$ are the same constants as in statements of Theorems~\ref{t:ui} and~\ref{t:vl}.

\section{Preliminary results}
\label{s:pre}
In this section, we collect some existing results that we exploit in the following section. In Subsection~\ref{sus:invman} we give some results on invariant manifolds for ordinary differential equations and in Subsections~\ref{sus:sss}, \ref{sus:rie} and \ref{sus:brie} we consider self-similar solutions of systems of conservation laws. More precisely, in Subsection~\ref{sus:sss} we investigate the structure of a general self-similar, distributional  solution having small total variation. In general, such a solution is not unique, so in Subsections~\ref{sus:rie} and~\ref{sus:brie} we describe how to construct a self-similar, distributional solution which has small total variation and which satisfies additional admissibility conditions. In particular,  in Subsection~\ref{sus:rie}, we deal with the Riemann problem~\eqref{S1: hyp},~\eqref{e:i:rie}, while in Subsection~\ref{sus:brie} we focus on the boundary Riemann problem~\eqref{S1: hyp},~\eqref{e:cdrie}.

\subsection{Invariant manifolds for ordinary differential equations}
\label{sus:invman}
In the following we rely on the notions of center, center-stable and uniformly stable manifold. For completeness, we present here the main properties that are exploited in the present paper. For an extended discussion, we refer the reader to the books by Katok and Hasselblatt~\cite{KHass} and by Perko~\cite{Perko} and to the lecture notes by Bressan~\cite{Bre:notes}. 

Let us consider the ordinary differential equation
\be
\label{e:ode}
     \frac{d V}{d y} = G(V), 
\ee     
where $V(y) \in \R^d$ for every $y$ and $G: \R^d \to \R^d$ is a smooth function which attains the zero value at some point. Without any loss of generality, we can assume that zero is an equilibrium point, namely $G(\vec 0) = \vec 0$. 

Let us assume that at least one eigenvalue of the Jacobian matrix $DG (\vec 0)$ has zero real part. Also, let us denote by $V^c$ the subspace of $\R^d$ spanned by vectors, the so-called \emph{generalized eigenvectors}, associated with these eigenvalues and let us denote by $n_c$ the dimension of $V^c$. The Center Manifold Theorem ensures that there exists a so-called center manifold $\mathcal M^c \subseteq \R^d$, which enjoys the following properties: first, $\mathcal M^c$ is locally invariant for~\eqref{e:ode}, namely, if $V_0 \in \R^d$, then the solution of the Cauchy problem 
$$
\left\{
\begin{array}{ll}
          d V/ d y = G(V)  \\
          V(0) = V_0
\end{array}
\right.
$$
satisfies $V(y) \in \mathcal M^c$ if $|y|$ is small enough. Also, there exists a small enough ball, $B_r (\vec 0)$, centered at $\vec 0$, such that, if $V(y) \in B_r (\vec 0)$ for every $y$, then $V(y) \in \mathcal M^c$ for every $y$. Finally, $\mathcal M^c$ has dimension $n_c$ and it is tangent to $V^c$ at $\vec 0$. 

Let us denote by $V^s$ the subspace of $\R^d$ spanned by the generalized eigenvectors of $DG(\vec 0)$ associated with the eigenvalues with strictly negative real part. The Center-Stable Manifold Theorem (see e.g.~\cite[Chapter 6]{KHass}) states that there exists a so-called center-stable manifold $\mathcal  M^{cs}$, tangent to $V^s \oplus V^c$ at the origin and locally invariant for~\eqref{e:ode}, which satisfies the following properties: there exists a sufficiently small ball $B_r (\vec 0)$ such that, if $V(y)$ is a solution of~\eqref{e:ode} satisfying $V(y) \in B_r (\vec 0)$ for every $y>0$, then $V(y) \in \mathcal M^{cs}$ for every $y$. 

Finally, let us consider the case that there exists a manifold $\mathcal E$ which consists of equilibria of~\eqref{e:ode}, namely $G(V) = \vec 0$ for every $V \in \mathcal E$. Let $E$ be the tangent space to $\mathcal E$ at the origin. Then there exists a manifold $\mathcal M^{us}_{\mathcal E}$, the so-called uniformly-stable manifold, which satisfies the following properties: first, $\mathcal M^{us}_{\mathcal E}$ is locally invariant for~\eqref{e:ode} and it is tangent to $V^s \oplus E$ at the origin. Also, there exists $B_r (\vec 0)$, a small enough ball centered at $\vec 0$, and a constant $c>0$ such that the following holds: if $V(y)$ is a solution of~\eqref{e:ode} satisfying 
$$
    \lim_{ y \to + \infty} | V(y) - \bar V| e^{c\,y /2 } = 0
$$
for some $\bar V \in B_r (\vec 0) \cap \mathcal E$, then $V(y) \in \mathcal M^{us}_{\mathcal E}$ for every $y$. The uniformly-stable manifold can be viewed as a particular example of \emph{slaving manifold}; see again~\cite{KHass} for a related discussion.  

\subsection{Self-similar solutions of systems of conservation laws}
\label{sus:sss}
Systems of conservation laws are invariant under uniform stretching of both space and time variables, and, thus, admit self-similar solutions in the form 
\be
	U(t,x)=V(\frac{x}{t})
\ee
where $V(\xi)$ is a measurable function of one variable. The above form implies that self-similar solutions are constant along straight-line rays emanating from the origin and, therefore, the study of such solutions is very important in the context of the Riemann problem, for which the initial data introduce a jump discontinuity at the origin. In this section, we consider self-similar solutions to Riemann problems and describe the structure of the solution. In particular, in this paper we consider initial-boundary value problems, so here we focus on the domain $(t, x) \in ]0, + \infty [ \times ]0, + \infty[$. Also, we assume that $V$ has small total variation, i.e.
\be 
\label{e:smbv}
        \mathrm{TotVar} V \leq C \delta, 
\ee
where $0 < \delta <<1$ and $C>0$ are as in the statement of Theorem~\ref{t:ui}. For a discussion on the main properties of functions in one space variable having bounded total variation, see Chapter 3.2 in the book by Ambrosio, Fusco and Pallara~\cite{AmbrosioFuscoPallara}.  

Having bounded total variation, $V$ admits at most countably many discontinuities, say located at the points $\{ j_m \}$, $m \in \N$. The left and right limits of $V$, as $\xi$ tends to a discontinuity point $j_m$, exist and are finite and we adopt the standard notation $V(j_m^+)$ and $V(j_m^-)$. Also, the limits as $\xi \to 0^+$ and $\xi \to + \infty$ exist and are finite and we denote them by 
\be
\label{e:lim}
        \lim_{\xi \to 0^+ } V(\xi) = \bar U \quad \text{and} \quad   \lim_{\xi \to + \infty } V(\xi) = U_0
\ee
respectively.  Following Dafermos~\cite[Section 9.1]{D}, we decompose $]0, + \infty[$ in disjoint sets as follows:
\be 
\label{e:deco} 
         ]0, + \infty[ =   \mathcal C \cup \mathcal W \cup \bigcup_{m \in \N} j_m,
\ee 
where $\mathcal C$ is the complement of the support of the measure $d V / d \xi$ and $\mathcal W$ is the (possibly empty) set of the points of continuity of $V$ that are not in $\mathcal C$. Since the measure $d V / d \xi$ vanishes on $\mathcal C$, then the function $V$ is constant on any interval included in $\mathcal C$. From~\eqref{e:smbv}, we also have 
\be 
\label{e:linfty}
        |V (\xi) - U_0 | \leq C \delta \quad \text{for every $\xi \in \, ]0, + \infty[$}. 
\ee
In the remaining part of this section, we focus on the boundary characteristic case~\eqref{e:c}-\eqref{e:lk}. The non characteristic case~\eqref{e:nchar} can be treated in an analogous way. Let $\lambda_i(V)$ be the $i$-th eigenvalue of the Jacobian matrix $DF(V)$ and by $r_i(V)$ a corresponding eigenvector with unit norm (the orientation is not important here). Here, we introduce some notation that we need in the following. We define
\be 
\label{e:pm}
\begin{split}
&         \lambda_i + : = \max_{  |V - U_0| \leq C \delta} \lambda_i (V) + M \delta, \qquad 
         \lambda_i-: = \min_{ |V - U_0| \leq C \delta} \lambda_i (V) - M \delta  
         \quad \text{for $i= k+1, \dots, n$}, \\
&       \qquad \qquad  \qquad \qquad \qquad \qquad   \lambda_k+ : =\max_{  |V - U_0| \leq C \delta} \lambda_k (V) + M \delta.
  \\
\end{split}
\ee
Here, the constant $M$ depends only on the flux $F$ and the initial data $U_0$ and its value is determined in the proof of Lemma~\ref{l:sfs}.
 
\begin{lemma}
\label{l:sfs}
         Let $U(t, x) = V (x /t)$ be a self-similar solution of the system of conservation laws 
         $$
             \partial_t U + \partial_x F(U) =0, \quad \text{ $(t, x) \in \, ] 0, + \infty [ \times ] 0, + \infty [ $}            
         $$
        and assume that $V$ satisfies~\eqref{e:smbv}. If the system is strictly hyperbolic and the constant $\delta$ in~\eqref{e:smbv} is small enough, then the following hold.
        \begin{enumerate}
        \item There are constant states 
        $U_k^\sharp$, \dots,$U_{n-1}^\sharp\in\R^n$, such that
\be\label{e:Vstates}
 V(\xi)=\left\{ \begin{array}{ll}
			U_i^\sharp &    \lambda_i+<\xi<\lambda_{i+1}-,\qquad i=k,\ldots,n- 1 \\
			U_0 &  \lambda_n+<\xi<+\infty
			\end{array}\right.
\ee
        \item At any point of discontinuity $j_m$, $m\in\N$, the Rankine-Hugoniot conditions hold, namely 
        \be
        \label{e:rh}
              F \big( (V(j^+_m) \big) - F \big( (j^-_m) \big) = 
              j_m \Big[ V (j^+_m) - V (j^-_m) \Big] \quad \text{for every $m \in \N$}.
        \ee
        \item Let $\mathcal W$ be the same as in~\eqref{e:deco}, then, for every $\xi \in \mathcal W$, 
        \be 
        \label{e:der}
                \lambda_i(V(\xi))=\xi 
         \ee       
         for some $i =1, \dots, n$. If $V$ is also differentiable at $\xi$, then 
         $d V / d \xi$ is parallel to $r_i \big( V (\xi) \big)$.  
          \item If $\xi$ is in the (possibly empty) inner part of $\mathcal W$ and 
         $$
             \nabla \lambda_i \big( V(\xi) \big) \cdot r_i \big( V(\xi) \big) \neq 0,
         $$
         then $V$ is differentiable at $\xi$.
          \end{enumerate}
\end{lemma}
\begin{proof}
Properties (2), (3) and (4) follow straightforwardly from the analysis in~\cite[Section 9.1]{D} and from Proposition 8.5 in Tzavaras~\cite{Tz}. Indeed, in the proof of that proposition one does not exploit that $V$ is obtained as a limit of a viscous approximation. Regarding the structure~\eqref{e:Vstates}, we assume that the constant $\delta$ is so small that the intervals
$ ]0, \lambda_k+[ $,  $ ]\lambda_{k+1}-, \lambda_{k+1}+[ , \dots  
    ]\lambda_{n}-, \lambda_n+[ $ are disjoint. One can first observe that from~\eqref{e:der}, it follows 
\be\label{Wsub}
   \mathcal W \subseteq \, ]0, \lambda_k+[ \,  \cup \, ]\lambda_{k+1}-, \lambda_{k+1}+[ \, \cup \, 
   \dots \cup \,  ]\lambda_{n}-, \lambda_n+[  \, .
\ee
Also, from the Rankine-Hugoniot conditions~\eqref{e:rh}, we get that 
$$
         \int_0^1 DF \Big( \theta V (j_m^+) + [1 - \theta]  V (j_m^-) \Big) d \theta \cdot   
          \Big[ V (j_m^+) -  V (j_m^-) \Big]  = j_m  \Big[ V (j_m^+) - V (j_m^-) \Big]
$$
and hence, $j_m$ is an eigenvalue of the matrix
$$
       \int_0^1 DF \Big( \theta V (j_m^+) + [1 - \theta]  V (j_m^-) \Big) d \theta.
$$
This implies that the constant $M$ in~\eqref{e:pm} can be chosen in such a way that
$$
     \bigcup_{m \in \N} j_m  \subseteq \, ]0, \lambda_k+[ \,  \cup \, ]\lambda_{k+1}-, \lambda_{k+1}+[ \, \cup \, 
   \dots \cup \,  ]\lambda_{n}-, \lambda_n+[  \, ,
$$
and, recalling~\eqref{Wsub}, Property (1) follows.
\end{proof}

\subsection{The solution of the Riemann problem}
\label{sus:rie}
Let us consider the Riemann problem for a system of conservation laws, namely
\be
\label{e:rie}
\left\{
\begin{array}{ll}
          \partial_t U + \partial_x F(U) =0\\
          U(0, x) =
          \left\{
          \begin{array}{ll}
                     U^+  & x > 0 \\
                     U^-   & x < 0\;,
          \end{array}
          \right.
\end{array}
\right.
\ee 
where $(t, x) \in [0, + \infty [ \times \R$ and $U(t, x) \in \R^n$. Here, $U^+$ and $U^-$ are two given constant states in $\R^n$. In \cite{lax}, Lax first introduced an admissibility condition on shocks, which is now named after him. Also, he constructed a distributional solution of \eqref{e:rie} satisfying this criterion. 
The analysis in \cite{lax} relied on the assumptions that $|U^+ - U^-|$  is small and that each characteristic field is either genuinely nonlinear or linearly degenerate. The main point is the construction of the so-called \emph{$i$-wave fan curve} through a given state $U^+$, which in the following is denoted by $T_i (s_i, U^+)$, $s_i$ being the variable parameterizing the curve, $i=1,\dots,n$. The curve $T_i( \cdot, U^+)$ attains value in $\R^n$ and satisfies the following property: given $s_i$ small enough, the Riemann problem 
$$
\left\{
\begin{array}{ll}
          \partial_t U + \partial_x F(U) =0  \\
          U(0, x) =
          \left\{
          \begin{array}{ll}
                      U^+  & x > 0 \\
                      T_i (s_i, U^+) & x < 0
          \end{array}
          \right.
\end{array}
\right.
$$
admits a self-similar solution   $V(x/t)= U(x, t)$, which is either a rarefaction, a single contact discontinuity or a single shock satisfying the Lax admissibility condition. Also, let $\lambda_i (U)$ denote the $i$-th eigenvalue of $DF (U)$, then there exists a small constant $\delta > 0$ such that function $V(\xi)$ is identically equal to $T_i (s_i, U^+)$ on the interval $\xi\in]- \infty, \lambda_i ( U^+) - \delta ]$ and is identically equal to $U^+$ on the interval $\xi\in[ \lambda_i ( U^+) + \delta, + \infty [ $. 

In the general case, when $U^- \neq T_i( s_i, U^+)$, a solution of~\eqref{e:rie} is obtained as follows: one defines the map 
\be 
\label{e:ciao}
    \chi(s_1, \dots ,s_n,  U^+ ) : = T_1 (s_1, T_2 (s_2, T_3 (\cdots ,T_n (s_n, U^+) \cdots )))  
\ee
and shows that it is locally invertible in a small enough neighborhood of $(s_1, \dots,s_n) = \vec 0$. Hence, given $U^-$, by imposing $U^- = \chi(s_1, \dots, s_n,  U^+ )$, one uniquely determines the 
value of $(s_1, \dots, s_n)$. The corresponding solution of~\eqref{e:rie} is made by a finite number of rarefaction waves, contact discontinuities  and shocks that satisfy the Lax admissibility condition. 

In~\cite{Liu:TAMS,Liu:rie}, Liu extended the analysis of Lax in~\cite{lax} to very general systems and introduced the so-called \emph{Liu admissibility condition}. Given $U^+ \in \R^n$, the so-called $i$-Hugoniot locus is a curve $W_i ( s, U^+)$, which is defined in a small enough neighborhood of $U^+$ and comprises all the states such that the Rankine-Hugoniot conditions 
$$
     F \big( W_i ( s, U^+) \big) - F(U^+) =  \sigma_i (s, U^+) \Big[   W_i ( s, U^+) - U^+ \Big]
$$
are satisfied for a certain real number $\sigma_i (s, U^+) $ close to $\lambda_i (U^+)$. Given a value $\bar s$, the shock between $W_i (\bar s, U^+)$ (on the left) and $U^+$ (on the right) is admissible in the sense of Liu if 
\begin{equation}
\label{e:liu}
          \sigma_i (\bar s, U^+) \leq \sigma (s, U^+) \quad \text{for every $s$ between $0$ and $\bar s$}. 
\ee
It should be mentioned that Liu admissibility condition can be viewed as an extension of Lax admissibility condition. See Dafermos~\cite[Chapter 8.4]{D} for a longer discussion.

In~\cite{BiaBrevv}, Bianchini and Bressan constructed the solution of~\eqref{e:rie} as the limit $\e \to 0^+$ of the vanishing viscosity approximation, namely the limit of solutions $U^\e$ to the parabolic system
\be
\label{e:rievv}
\left\{
\begin{array}{ll}
          \partial_t \ue + \partial_x F(\ue) = \e \partial_{xx} \ue \\
          \ue (0, x) =
          \left\{
          \begin{array}{ll}
                     U^+  & x > 0 \\
                     U^-   & x < 0\;.
          \end{array}
          \right.
\end{array}
\right.
\ee 
Since in this work we exploit ideas and techniques from~\cite{BiaBrevv}, we now go over the construction of the $i$-wave fan curve therein. If $s_i <0$, then one considers the following fixed point problem, defined on the interval $[s_i, 0]$:
\begin{equation}
\label{e:biabre}
         \left\{
         \begin{array}{lll}
                   \displaystyle{U_i (\tau) = U^+ + \int_0^{\tau}
                    \tilde{r}_i \big( U_i(z), v(z), \sigma_i (z)\big) dz} \\
                   \displaystyle{ v(\tau) = f (\tau) - 
                   \mathrm{conv} f (\tau)} \\
                   \displaystyle{ \sigma_i (\tau)  = 
                   \frac{d}{d \tau}  \mathrm{conv} f (\tau)} . \\
          \end{array}
         \right.
\end{equation}
In the above problem, $\tilde r_i$ is the so-called \emph{generalized eigenvector} defined in~\cite[Section 4]{BiaBrevv} by relying on the Center Manifold Theorem. Also, $f$ is the \emph{generalized flux} defined by
$$
      f (\tau) = \int_0^{\tau} \tilde{\lambda}_i \big( U_i(z), v(z),  \sigma_i(z) \big) dz, 
$$
where $\tilde{\lambda}_i $ is the \emph{generalized eigenvalue}
$$
  \tilde{\lambda}_i ( U_i, v,   \sigma_i ) = DF(U) \tilde r_i ( U_i, v,   \sigma_i )  \cdot \tilde r_i ( U_i, v,   \sigma_i )\;.
$$
Finally, $\mathrm{conv} f$ represents the convex envelope of $f$, given by
\begin{equation}
\label{e:conv}
    \mathrm{conv} f (\tau) = \sup \big\{g(\tau): \;
     \text{$g$ convex and $g(z) \leq f(z) \; \forall \, z \in [s_i, 0]$}  \big\} .
\end{equation}
One can show that the fixed point problem~\eqref{e:biabre} admits a unique solution $\big(U_i(\tau) , v (\tau) , \sigma_i (\tau) \big)$ and hence the $i$-th wave fan curve $T_i(\cdot, U^+)$ is defined by setting $T_i(s_i, U^+): = U(s_i)$. If $s_i >0$, the construction is entirely similar, the only difference being that one has to take in~\eqref{e:biabre} the concave envelope of $f$ instead of the convex one. 

The basic idea in the construction of the $i$-wave fan curve in~\cite{BiaBrevv} is the following: let us fix $s_i <0$ and consider the solution of the fixed point problem~\eqref{e:biabre}. Assume that $v(\tau) >0$ on a given subinterval $]\alpha, \beta [ \subseteq [s_i, 0]$, and that $v (\alpha) = v(\beta) =0$. Then, one can show that a solution of the Riemann problem
$$
\left\{
\begin{array}{ll}
          \partial_t U + \partial_x F(U) =0 \\
                    U(0, x) =
          \left\{
          \begin{array}{ll}
                     U (\beta)  & x > 0 \\
                     U (\alpha)    & x < 0
          \end{array}
          \right.
\end{array}
\right.
$$
is a shock (or a contact discontinuity) which is also admissible in the sense of Liu. Conversely, assume that 
$v(\tau) = 0$ on a given subinterval $]a, b [ \subseteq [s_i, 0]$, then one can show that a solution of the Riemann problem
$$
\left\{
\begin{array}{ll}
          \partial_t U + \partial_x F(U) =0 \\
                    U(0, x) =
          \left\{
          \begin{array}{ll}
                     U (b)  & x > 0 \\
                     U (a)    & x < 0
          \end{array}
          \right.
\end{array}
\right.
$$
is either a rarefaction wave or a contact discontinuity which is admissible in the sense of Liu. The solution is a contact discontinuity if $\nabla \lambda_i (U) \cdot r_i (U) \equiv 0$ on $ ]a, b[$. In general, the state $U(s_i)$ is connected to $U^+$ by a sequence of rarefaction waves and of shocks (or contact discontinuities) satisfying the Liu admissibility condition. 

Once the wave fan curves $T_1, \dots, T_n$ are defined, the argument in~\cite{BiaBrevv} works as in Lax~\cite{lax}: one defines the function $\chi$ as in~\eqref{e:ciao} and shows that it is locally invertible, hence by setting $U^- = \chi (s_1, \dots ,s_n, U^+)$ one uniquely determines 
$(s_1, \dots ,s_n)$, provided $|U^+ - U^-|$ is small enough. The limit $\e \to 0^+$ of~\eqref{e:rievv} is then given by 
$$
    U(t, x) =
    \left\{
    \begin{array}{ll}
               U^- & x \leq \sigma_1 (s_1)\,t\\
               U_i(s_i)   & \sigma_i (0) t\le x \le \sigma_{i+1}(s_{i+1})\,t \qquad i=1,\dots,n-1 \\
               U_i (\tau) & x = \sigma_i (\tau) t\qquad i=1,\dots,n  \\
               U^+ & x \ge \sigma_n (0)\, t .
    \end{array}
    \right.
$$
In the previous expression, $U_i$ and $\sigma_i$ are obtained by considering the solution of the fixed point problem~\eqref{e:biabre}. It should be noted that in the case that each characteristic field is either genuine nonlinear or linearly degenerate, the wave fan curve constructed by Bianchini and Bressan in~\cite{BiaBrevv} is the same as in Lax~\cite{lax}.

In~\cite{Bia:riemann}, Bianchini extended the construction in~\cite{BiaBrevv} by considering more general approximations of the Riemann problem~\eqref{e:rie}. In particular, he handled the viscous approximation   
\be
\label{e:b}
     \partial_t \ue + \partial_x \Big[ F(\ue) \Big] = \partial_x \Big[ B(\ue) \partial_x \ue \Big],
\ee
under some hypotheses on the matrix $B(U)$. He also established a uniqueness result for~\eqref{e:rie}, namely that there exists a unique solution $U$ of~\eqref{e:rie} which is obtained by patching together a sequence of rarefaction waves, shocks and contact discontinuities in such a way that the corresponding speed is increasing as one moves from the left to the right and in such a way that all shocks and contact discontinuities are admissible in the sense of Liu. In particular, this implies that the solution of~\eqref{e:rie} obtained by taking the limit of the viscous approximation~\eqref{e:b} does not change if one changes the viscosity matrix $B(U)$. 

Also, Dafermos~\cite[Chapter 9]{D} provided a characterization of the limit $\e \to 0^+$ of
the self-similar viscous approximation 
$$
\left\{
\begin{array}{ll}
           \partial_t \ue + \partial_x \Big[ F(\ue) \Big] =  \e t \partial_{xx} \ue \\
          \ue(0, x) =
          \left\{
          \begin{array}{ll}
                     U^+  & x > 0 \\
                     U^-   & x < 0
          \end{array}
          \right.
\end{array}
\right.
$$
using a center manifold technique. In particular, he constructed the $i$-th wave fan curve for the self-similar approximation and showed that it is the same as the $i$-th wave fan curve $T_i$ for the viscous approximation~\eqref{e:b}. See also the analysis in Tzavaras~\cite{Tzavaras:JDE, Tz} and Joseph and LeFloch~\cite{JoLeF:nc}.

\subsection{The solution of the boundary Riemann problem}
\label{sus:brie}
Consider the viscous approximation of a boundary Riemann problem
\be
\label{e:brie}
\left\{
\begin{array}{ll}
           \partial_t \ue + \partial_x \big[ F(\ue) \big] = \e \partial_x \big[ B(\ue) \partial_x \ue \big] \\
          \ue(0, x) =U_0, \quad x>0 \\
          \ue(t, 0) = U_D , \quad t>0,
\end{array}
\right.
\ee 
where $U_0$ and $U_D$ are two given constant states in $\R^n$. 
For now, we focus on the case when the viscosity matrix $B(U)$ is invertible for every $U$. 

The analysis of the boundary Riemann problem was initiated in Dubois and LeFloch~\cite{DubLeF}. One of the main difficulties in studying the limit of~\eqref{e:brie} is the presence of \emph{boundary layer} phenomena. Convergence results have been proved in some special cases, see e.g. Gisclon \cite{Gis} and Ancona and Bianchini~\cite{AnBia}. 

The way boundary layers come into play is the following. Assume that $\ue $ converges as $\e\to0^+$ (in a suitable topology) to a limit function $U$ having bounded total variation. One does not expect that $U$ satisfies the boundary condition in~\eqref{e:brie}, namely, in general, $\lim_{x \to 0^+} U(t, x) \neq U_D$. Note that this limit is well defined since the solution $U$ is assumed to be of bounded total variation with respect to the $x$ variable. Let $\bar U$ denote the trace $$\bar U = \lim_{x \to 0^+} U(t, x),$$ then one expects that the relation between $\bar U$ and $U_D$ is the following: there exists a \emph{boundary layer} $W$ satisfying 
\be 
\label{e:bl1}
\left\{
\begin{array}{ll}
          B(W) W' = F(W) - F( \bar U) \\
          W(0) = U_D \quad \text{and} \quad \lim_{y \to + \infty} W(y) = \bar U.  
\end{array}
\right.
\ee 
Gisclon~\cite{Gis} and Gisclon and Serre~\cite{GisSerre} established, among other results, a detailed analysis of the boundary layers for a general viscous approximation when the viscosity matrix $B$ is invertible. See also Joseph and LeFloch~\cite{JosephLeFloch:ARMA} for the analysis of the boundary layers coming from other approximations. 

Concerning the analysis of the limit of~\eqref{e:brie}, the case of identity matrix, $B(U) \equiv I$, was considered by Ancona and Bianchini~\cite{AnBia}, who established convergence results and provided a characterization of the limit. The analysis of the limit was then extended in Bianchini and Spinolo~\cite{BiaSpi} to the case of more general viscosity matrices, including the case that $B(U)$ is singular, which is natural in view of most of the physical applications. For the convenience of the reader, we present the analysis in~\cite{BiaSpi}, since we use it in Section~\ref{s:uni}. In the case that the viscosity matrix $B$ is singular, the initial-boundary value problem~\eqref{e:brie} may be ill-posed and hence one has to use a more complex formulation of the boundary condition. Here, for simplicity, we discuss only the case of an invertible viscosity matrix. First, we focus on the case of a non-characteristic boundary, namely we assume that all the eigenvalues of the Jacobian matrix $DF(U)$ are bounded away from $0$. The boundary characteristic case is considered afterwards.

Let $\bar U$ denote as before the hyperbolic trace, then we expect that $\bar U$ and the Cauchy datum $U_0$ are connected by a sequence of rarefaction waves, shocks, and contact discontinuities having \emph{positive} speed and being admissible in the sense of Liu. Namely,
\be 
\label{e:trace}
    \bar U = T_{n-p+1} (
    s_{n-p+1}, T_{n-p+2}(
    s_{n-p+2}, T_{n-p+3} (\cdots , T_n(s_n, U_0 ) \cdots )
    )
    )
    \ee
for some values $s_{n-p+1}, \dots ,s_n$.
In the previous expression, $p$ denotes the number of positive eigenvalues of $DF(U)$, while $T_i(s_i, U_0)$ is the same $i$-wave fan curve as in Bianchini~\cite{Bia:riemann}. We recall that the construction of $T_i$ is presented in Section~\ref{sus:rie}. To express the relation between $\bar U$ and $U_D$, we exploit~\eqref{e:bl1} and the Stable Manifold Theorem. This theorem ensures that there exists a manifold, the so-called \emph{stable manifold}, having dimension equal to the number of eigenvalues of $B^{-1} (\bar U) DF (\bar U)$ with negative real part. Under some assumptions on $B$ (see~Lemma~\ref{l:same} in Section~\ref{sss:hyp}), this number is exactly $n-p$. Moreover, all the functions $W(y)$  
satisfying $B(W) W' = F(W) - F(\bar U)$ and $\lim_{y \to + \infty} W(y) = \bar U$ lie on the stable manifold. Hence, if be denote by $\varphi(s_1, \dots, s_{n-p}, \bar U)$ a function parameterizing the manifold, by combining the above argument with~\eqref{e:trace}, we get that 
\be
\label{e:ud}
     U_D = \varphi \Big( s_1, \dots, s_{n-p}, T_{n-p+1} (
    \cdots , T_n(s_n, U_0 ) \cdots )
    )
    )
    \Big).
\ee
Under some assumptions on the matrix $B$, the function
$$
   (s_1, \dots, s_n) \mapsto \varphi \Big( s_1, \dots, s_{n-p}, T_{n-p+1} (\cdots , T_n(s_n, U_0 ) \cdots ) 
    \Big)
$$
is locally invertible and hence, by imposing~\eqref{e:ud}, we determine uniquely the value of $(s_1, \dots, s_n)$. By arguing as in Section~\ref{sus:rie}, the value $U(t,x)$ of the limit of~\eqref{e:brie} can be determined for almost every $(t,x)$.  
  
Let us now consider the limit $U$ of~\eqref{e:brie} in the boundary characteristic case, namely we assume that the $k$-th eigenvalue of the Jacobian matrix $DF(U)$ can attain the value $0$, i.e. $\lambda_k (U) \sim 0$. The analysis in this case is much more delicate because there is not strict separation between the boundary layers on one side and the rarefaction waves, the shocks and the contact discontinuities with positive speed on the other side. 

For the boundary characteristic case, the construction in~\cite{BiaSpi} works as follows. There exists a state, denoted by $U_k^{\sharp}$, which is connected to $U_0$ by rarefaction waves, shocks and contact discontinuities of the families $k+1, \dots, n$:
\be
\label{e:usharp}
    U_k^{\sharp}= T_{k +1} (
    s_{k+1 }, T_{k+ 2}(
    \cdots , T_n(s_n, U_0 ) \cdots )
    ),
\ee
for some $s_{k+1},\dots,s_n$. To complete the analysis, instead of the $k$-wave fan curve $T_k$, one employs the \emph{characteristic wave fan curve} $\tchar ( s_k, U_k^{\sharp})$, which is constructed as follows. If $s_k<0$, we consider the fixed point problem  
\begin{equation}
\label{e:char1}
         \left\{
         \begin{array}{lll}
                   \displaystyle{U_k (\tau) = U^{\sharp} + \int_0^{\tau}
                    \tilde{r}_k \big( U_k(z), v_k(z), \sigma (z)\big) dz} \\
                   \displaystyle{ v_k(\tau) = f (\tau) - 
                   \mathrm{monconv} f (\tau)} \\
                   \displaystyle{ \sigma (\tau)  = \frac{1}{d} 
                   \frac{d}{d \tau}  \mathrm{monconv} f (\tau)} , \\
          \end{array}
         \right.
\end{equation}
where $\tau\in[s_k, 0]$. Here, $\tilde{r}_k$ is the same as in~\eqref{e:biabre} and $f$ denotes again a ``generalized flux". Finally, $\mathrm{monconv} f$ represents the monotone convex envelope of $f$, namely
\begin{equation}
\label{e:monconv}
    \mathrm{monconv} f (\tau) = \sup \big\{g(\tau): \;
     \text{$g$ convex, non decreasing, $g(z) \leq f(z) \; \forall \, z \in [s_k, 0]$}  \big\} .
\end{equation}
By relying on the Contraction Map Theorem, one can show that there exists a unique continuous solution $\big( U_k(\tau), v_k(\tau), \sigma (\tau) \big)$ of~\eqref{e:char1} which is confined in a small enough neighborhood of 
$\big( U_k^{\sharp}, 0, \lambda_k (U_k^{\sharp} ) \big)$.  One eventually defines the characteristic wave fan curve by setting $ \tchar (s_k, U_k^{\sharp}) : = U(s_k)$. If $s_k>0$, the construction is entirely similar, the only difference being that in~\eqref{e:char1} one has to take the monotone concave envelope of $f$ instead of the monotone convex one. 

The idea behind the construction sketched above is the following. Let us fix $s_k<0$ and consider the solution of~\eqref{e:char1}, defined on the interval $[s_k, 0]$ and set 
$$
   \bar s : = \max \{ \tau \in [s_k, 0] \; :\; \sigma (\tau) =0 \}
$$ 
and 
$$\underline s: = \sup \{ \tau \in [s_k, 0] \; \textrm{s.t. $\sigma (\tau) =0$ and $v_k (\tau)>0$} \}.
$$ 
Clearly, $\underline s\le \bar s$.
Let us focus on the most interesting case by assuming $s_k < \underline s < \bar s < 0$. Then $U(\bar s)$ (on the left) is connected to $U(0)= U_k^{\sharp}$ (on the right) by a sequence of rarefaction waves and shocks (or contact discontinuities) satisfying the Liu admissibility condition and having strictly positive speed. Also, $U(\underline s)$ (on the left) is connected to $U(\bar s)$ (on the right) by a sequence of contact discontinuities that are admissible in the sense of Liu and have speed $0$. Moreover, by construction, $v_k(\tau) > 0 $ and $\sigma (\tau) =0$ for $\tau\in[s_k, \underline s[$. Hence, the Cauchy problem
\be
\label{e:chv}
   \left\{
    \begin{array}{ll}
               \displaystyle{\frac{d y }{d \tau} = \frac{1}{v_k (\tau)}} \\ 
               y (s) =0
    \end{array}       
   \right.
\ee
defines a change of variables: $\tau\mapsto y$, $[s_k, \underline s [ \to [ 0, + \infty [$. One can show that  $W_k (y) : = U_k (\tau (y))$ is a solution of the equation $B(W) W' = F(W) - F ( U_k(\underline s))$ satisfying $W_k(0) = \tchar (s_k, U_k^{\sharp})$ and $\lim_{y \to + \infty} W_k (y) = U_k (\underline s)$. 

To complete the construction, one defines a manifold, transversal to the vector  $\tchar (s_k, U_k^{\sharp}) - U_k^{\sharp}$, having dimension $k-1$ and satisfying the following property: if $W_s (0)$ lies on this manifold, namely $W_s(0)=\omega(s_1,\dots,s_{k-1})$ for some parameterization $\omega$ of the manifold, then the solution $W_{s}$ of the equation 
\be\label{e:ws}
               \Big[ B(W_{s}+ U_k^{\sharp} ) [W_{s} + U_k^{\sharp} ]' \Big]' = 
               \big[ F(W_{s}+ U_k^{\sharp}) \big]'
               \ee
decays exponentially fast to the equilibrium point $\vec 0$ as $y \to +\infty$. Moreover, for fixed $(s_1,\dots,s_{k-1})$ and $s_k$, one singles out a function $W_p$, such that by setting 
\be
\label{e:w}
     W(y) = W_s(y) + W_k (y)+W_p(y)\,,
\ee
one obtains a solution $W$ of system 
$ [ B(W ) W ' ]' =               [ F(W) ]'$ satisfying $\lim_{y \to + \infty} W(y) = U_k( \underline s)$, where $U_k$ is the solution of~\eqref{e:char1}. In~\eqref{e:w}, $W_k(y) $ is the above-mentioned function obtained by considering the change of variables~\eqref{e:chv} and $W_s(y)$ is a solution of~\eqref{e:ws} satisfying $W_s(0)=\omega(s_1,\dots,s_{k-1})$. One can determine $W_p$ 
in such a way that $W_p$ is small (in the sense specified in~\cite{BiaSpi}) with respect to $W_s$ and $W_k$, namely it can be regarded as a perturbation term.  By construction, $W(y)$, and hence the value $W(0)$, depends on the vector $(s_1, \dots, s_{k-1})$, on the number $s_k$ and on the state $U_k^{\sharp}$. One sets 
\be
\label{e:psi}
W(0): = \psi (s_1, \dots s_k, U_k^{\sharp})
\ee and then, by recalling~\eqref{e:usharp}, one defines the map
$$
   (s_1, \dots, s_n, U_0)  \mapsto
   \psi \Big( 
             s_1, \dots, s_{k-1}, s_k, T_{k-1} \big( 
                          \cdots,  T_n(s_n, U_0)  \big) 
              \Big), \\
$$
which, under suitable assumptions on the matrix $B$, is locally invertible. Hence, by imposing 
$$
U_D= \psi(s_1, \dots s_n, U_0),$$ 
one uniquely determines $(s_1, \dots ,s_n)$ and hence the limit of~\eqref{e:brie}. This completes the construction in~\cite{BiaSpi} for the boundary characteristic case.

Finally, we underline that in~\cite{JosephLeFloch, JoLeF:nc}, Joseph and LeFloch had previously constructed a solution of 
the boundary Riemann problem~\eqref{S1: hyp},~\eqref{e:cdrie}, by taking the limit $\e \to 0^+$ of the self-similar viscous approximation 
$$
\left\{
\begin{array}{ll}
           \partial_t \ue + \partial_x \big[ F(\ue) \big] = \e t \partial_x \big[ B(\ue) \partial_x \ue \big] \\
          U(0, x) =U_0 \quad \text{and} \quad U(t, 0) = U_D ,
\end{array}
\right.
$$ 
and they established the convergence under suitable hypotheses on $B$. Also, they provided a precise description of the limit, in both the non-characteristic and the characteristic case. In the boundary characteristic case, most of the analysis relied on the assumption that the $k$-characteristic field is genuinely nonlinear.

\section{Proof of Theorem~\ref{t:ui}}
\label{s:uni}
The existence part of the theorem is proved in~\cite[Section 3]{BiaSpi}. More precisely, by applying Theorems $3.1$ and $3.2$ in~\cite{BiaSpi}, we get the following. There exist constants $C$ and $\delta_1$, $\delta_1$ sufficiently small, such that if
\be\label{e:3.1a} |U_0-U_D|\le \delta_1,\ee
then there exist a function $V$ satisfying conditions $(1)$--$(5)$ in the statement of Theorem~\ref{t:ui} with $\delta=\delta_1$. Here, we focus on the issue of uniqueness. We distinguish between two cases: if the boundary is characteristic, namely if one eigenvalue of the Jacobian matrix $DF(U)$ can attain the value $0$, then the proof relies on two preliminary results, Lemmas~\ref{l:dec} and~\ref{t:u}, which are introduced in Sections~\ref{s:dec} and~\ref{sus:uwave} respectively. 
In the case when the boundary is non characteristic, the proof of Theorem~\ref{t:ui} is actually much easier and it follows from an argument in Joseph and LeFloch~\cite{JosephLeFloch}. Both cases are discussed in Section~\ref{sus:conchar}. 

\subsection{A decomposition result}
\label{s:dec}
In this subsection, the principal result is a decomposition of the boundary layers and its construction is given in Lemma~\ref{l:dec}. 

Throughout this subsection and the following one, we focus on the boundary characteristic case, i.e. we assume that
\be
\label{e:c1}
    \lambda_1 (U) < \dots < \lambda_{k-1}(U) < -c <0 < c < \lambda_{k+1}(U) < \dots < \lambda_n (U)   
\ee
and
\begin{equation}
\label{e:lk1}
   | \lambda_k (U) | \leq K \delta_1\;,
\end{equation}
for some constants $K>0$ and the same constant $\delta_1$ as in~\eqref{e:3.1a}. 

To begin with, we introduce some notations: we set 
$$
     A(U, U') = DF (U) - \Big[ B(U) \Big]'
$$
and consider boundary layers, namely functions $W$ satisfying 
\begin{equation}
\label{e:bl}
    \left\{
    \begin{array}{ll}
    W'' = B^{-1} (W) A( W, W') W' \\
    \lim_{y \to + \infty} W(y) \; \textrm{exists and it is finite}. \\
   \end{array}
   \right. 
\end{equation}
The above ordinary differential equation can be rewritten as a first order system    
\begin{equation}
\label{e:sy}
    \left\{
    \begin{array}{lll}
    W'  = p \\
    p'   = B^{-1} (W) \big[ A( W, p)  p - \sigma p \big] \\
    \sigma' = 0.
   \end{array}
   \right. 
\end{equation}
The reasons why we introduce the ``fake unknown" $\sigma$ would become clear in Sections~\ref{sus:uwave} and~\ref{sus:conchar}. 

For any given state $U^{\sharp} \in \R^n$, we can linearize~\eqref{e:sy} about the equilibrium point $\big( U^{\sharp}, \vec 0, \lambda_k (U^{\sharp}) \big)$. By applying the Center Manifold Theorem, we infer that there exists a center manifold $\mathcal M^c$ having dimension $n + 2$. For some more details of the construction, we refer the reader to {\sc Step 1} in the proof of the lemma below. Note that the center manifold is not unique, but we can arbitrarily fix one.

\begin{lemma} 
\label{l:dec} 
         Let conditions~\eqref{e:eig}--\eqref{e:diss} in Section~\ref{sss:hyp} be satisfied and let $U^{\sharp}$ be a state in $\R^n$ satisfying $|U^\sharp-U_0|\le C\delta_1$. Then there exist constants $C_1$, $C_2$ and $\delta$, $\delta\le\delta_1$, such that the following holds. Assume that ~\eqref{e:c1} and~\eqref{e:lk1} hold, and that $W$ is a solution of~\eqref{e:bl} satisfying
\begin{equation}\left\{
\label{e:bd}
   \begin{array}{lll}
   |W(y) - U^{\sharp} | \leq C  \delta,\qquad \text{for every}\,\,y\in[0,+\infty[ \\
   \lim_{y \to + \infty} W'(y) =\vec 0 \\ 
   |W'(y)| \leq C\delta \qquad  \forall \; y \in [0, + \infty[. 
   \end{array}
   \right.
\end{equation}
         Then $W$ can be decomposed as 
          \begin{equation}
          \label{e:dec}
              W (y) = U_k (y) + U_s (y) + U_p (y) \quad \text{for every $y \in [0, + \infty [ $},  
          \end{equation}
          where 
          \begin{enumerate}
          \item $U_k(y)$ lies on the center manifold $\mathcal M^c$ for every $y\in[0,+\infty[$ and satisfies $\lim_{y \to + \infty} U'_k (y) = \vec 0$.
          \item $U_s$ satisfies 
          \begin{equation}
          \label{e:uesse}
              |U_s (y) | \leq C_1 \delta e^{-c y  /  2} \; \textrm{for every $y \in [0, + \infty[$}
          \end{equation}
          and $B (U_s + U^{\sharp}) U''_s (y) =  A (U_s + U^{\sharp}, U'_s) U'_s $.
          \item the perturbation term $U_p$ is small with respect to the previous ones, more precisely
          \begin{equation}
          \label{e:upi}
              |U_p (y)| \leq C_2 \delta^2 e^{-c y  /  4} \; \textrm{for every $y \in [0, + \infty[$}.
          \end{equation}
          \end{enumerate}
          The constant $c$ in~\eqref{e:uesse} and~\eqref{e:upi} is the same as in~\eqref{e:c1}. 
\end{lemma}  
\begin{proof}
We exploit the construction in~\cite{AnBia} and~\cite{BiaSpi} and we rely on the notions of center, center-stable and uniformly stable manifold. For an overview, see Section~\ref{sus:invman} and the references therein. The proof is established in several steps.

\noindent
{\sc Step 1} We linearize system~\eqref{e:sy} about the equilibrium point $\big( U^{\sharp}, \vec 0, \lambda_k (U^{\sharp}) \big)$ and we obtain the $(2n+1) \times (2n+1)$ matrix
$$
    \left(
    \begin{array}{ccc}
              \mathbf{0} & {I_n}  & \mathbf{0} \\
              \mathbf{0}  & B^{-1} (U^{\sharp}) \big[ A( U^\sharp, \vec 0)   - \lambda_k (U^{\sharp})
               I_n \big] & \mathbf{0}  \\
               \vec 0^t & \vec 0^t & 0 \\
    \end{array}
    \right),
$$
where $\vec 0^t$ is the zero row vector in $\R^n$ and $\mathbf{0}$ and $I_n$ denote respectively the null and the identity matrices with dimension $n \times n$. 
By relying on Lemma~\ref{l:same}, we get that the matrix 
$$B^{-1} ( U^\sharp) \Big[ A (U^\sharp, \vec 0) -\lambda_k (U^\sharp) I_n \Big]$$ 
has the same number of eigenvalues with non-positive real part as the matrix $ \big( A (U^\sharp, \vec 0) -\lambda_k (U^\sharp) I_n \big)$. Consider the center-stable space, that is, the subspace of $\R^{2n +1}$ generated by the generalized eigenvectors corresponding to eigenvalues with non positive real part. Because of~\eqref{e:c1} and~\eqref{e:lk1}, the center-stable space has dimension $n+ 1+ k$ and therefore, by the Center-Stable Theorem, there exists a $(n+ 1+ k)$-dimensional center-stable manifold $\mathcal M^{cs}$ containing any solution of~\eqref{e:sy} that, for $y \in [0, + \infty[$, is confined in a small enough neighborhood of $\big( U^{\sharp}, \vec 0, \lambda_k (U^{\sharp}) \big)$. Let us fix a function $W$ satisfying~\eqref{e:bl} and~\eqref{e:bd}: if $\delta$ is small enough, then $\big( W(y), p(y)= W'(y), 0)$ lies on $\mathcal M^{cs}$. \\
{\sc Step 2} In the following, we use the basis composed by the eigenvectors $r_1, \dots, r_k$  of $A(U^{\sharp}, \vec 0)$ in $\R^n$. It should be noted that these eigenvectors are linearly independent because, by assumption, $A$ has $n$ real and distinct eigenvalues. By arguing as in~\cite[pages 45-47]{BiaSpi}, we get that $W'(y)$ admits the representation 
\begin{equation}
\label{e:uprime}
   W'(y) = R_{cs} \Big( W(y), V_{cs}(y), 0 \Big) V_{cs}(y),
\end{equation}
where $V_{cs}(y)$ is a function taking values in $\R^k$ and $R_{cs}$ is a function taking values in the space of $n \times k$ matrices and satisfying the following: $R_{cs} \big( U^{\sharp}, \vec 0, \lambda_k (U^{\sharp}) \big)$ is the matrix whose columns are $r_1 \cdots r_k$. \\
{\sc Step 3} We now employ the same ``diagonalization procedure" as in~\cite[pages 47-51]{BiaSpi} and we write 
\begin{equation}
\label{e:vcs}
     V_{cs} (y) = R_s \Big( W(y), V_s(y), 0 \Big) V_s (y) + \check{r}_k \Big( W(y), v_k (y), 0 \Big) v_k (y),
\end{equation}
where $V_s$ and $v_k$ take values in $\R^{k-1}$ and $\R$, respectively . The functions $\check{r}_k$ and $R_s$ take values in $\R^k$ and in the space of $k \times (k-1)$ matrices respectively. Also, from the construction in~\cite{BiaSpi} it follows that, for every $W$ and $v_k$, the point
$$
    \Big( W, R_{cs} (W, \vec 0, v_k, 0) \check r_k (W, v_k, 0)v_k , 0 \Big)  
$$
lies on the center manifold $\mathcal M^c$ defined just before the statement of Lemma~\ref{l:dec}. Moreover, for every $W$ and $V_s$, the point
$$
    \Big( W, R_{cs} (W, V_s, 0, 0)R_s (W, V_s, 0) V_s,  0 \Big)  
$$
belongs to the \emph{uniformly stable manifold} $\mathcal M^{us}_{\mathcal E}$ of system~\eqref{e:sy} with respect to the manifold of equilibria 
$$
    \mathcal E = \left\{ (W, \vec 0, \sigma), \; W \in \R^n ,\,\sigma\in\R\right\}. 
$$
We refer the reader to Section~\ref{sus:invman} for the notion of uniformly stable manifold. Finally, from the construction in~\cite{BiaSpi} follows that $\check{r}_k \big( U^{\sharp}, 0, \lambda_k (U^{\sharp}) \big) = r_k$ and $R_s \big( U^{\sharp}, \vec 0, \lambda_k (U^{\sharp}) \big)$ is the matrix whose columns are $r_1 \cdots r_{k-1}$. Hence, if $W(y)$ is fixed, the map $(V_s, v_k) \mapsto V_{cs}$ is locally invertible. In view of the above analysis, we consider $W(y)$ to be the solution of~\eqref{e:bl} and we first obtain the function $V_{cs}(y)$ from~\eqref{e:uprime} and then, using the local invertibility of the map, we get the functions $V_s(y)$ and $v_k(y)$ satisfying~\eqref{e:vcs}.  \\
{\sc Step 4.} By relying again on the ``diagonalization procedure" in~\cite[pages 47-51]{BiaSpi}, we get that $\big( W (y), v_k(y), V_s(y))$ satisfy 
\begin{equation}
\label{e:cs}
         \left\{
         \begin{array}{llll}
                  W'= R_{cs} (W, v_k \check{r}_k + R_s V_s, 0) 
                  \check{r}_k v_ k + R_{cs} (W, v_k \check{r}_k +
                   R_s V_s, 0)  R_s V_s \\
                  v_k'= \phi_k ( W, v_k, V_s, 0) v_k \\
                   V_s'= \Lambda ( W, v_k, V_s, 0 ) V_s, \\
         \end{array}
         \right.
\end{equation}  
where $\phi_k$ is a real-valued function satisfying $\phi_k ( U^\sharp, 0, \vec 0, \lambda_k(U^\sharp))=0 $, while $\Lambda$ attains values in the space of $(k-1) \times (k-1)$ matrices and 
$\Lambda (U^{\sharp}, 0, \vec 0, \lambda_k (U^\sharp))$ is a diagonal matrix having entries $\lambda_1 (U^{\sharp}) \dots \lambda_{k-1} (U^{\sharp})$. 
\\
{\sc Step 5.} We now illustrate our strategy to complete the proof of the lemma, the details are provided in the following step. 

We write 
\begin{equation}
\label{e:v}
        V_s (y) = \tilde V_s (y)+ V_p (y) \quad  \text{and} \quad 
        v_k (y) = \tilde v_k (y) + v_p (y), 
\end{equation}        
where $(\tilde V_s, \tilde v_k)$ are the ``principal components" and $(V_p, v_p)$ the ``perturbations". We impose that the components 
$U_s$ and $U_k$ in~\eqref{e:dec} satisfy 
\begin{equation}
\label{e:k}
     \left\{
         \begin{array}{llll}
                  U'_k= R_{cs} (U_k, \tilde v_k \check{r}_k, 0) 
                  \check{r}_k (U_k, \tilde v_k,  0)  \tilde v_ k  \\
                  \tilde v_k'= \phi_k ( U_k, \tilde v_k, \vec 0, 0) \tilde v_k \\
         \end{array}
         \right.
\end{equation}
and 
\begin{equation}
\label{e:s}
 \left\{
         \begin{array}{llll}
                  U_s'= R_{cs} (U_s + U^{\sharp},  R_s \tilde V_s, 0) 
                   R_s  (U_s + U^{\sharp}, \tilde V_s,  0)  \tilde  V_s \\
                   \tilde V_s'= \Lambda ( U_s + U^{\sharp}, 0, \tilde V_s, 0 ) \tilde V_s \\
         \end{array}
         \right.
\end{equation}
respectively. 

Hence, by exploiting~\eqref{e:dec} and~\eqref{e:v} and by plugging~\eqref{e:k} and~\eqref{e:s} in~\eqref{e:cs} we obtain the equation for $(U_p, v_p, V_p)$. Then, by imposing that both $U_p$ and $U_s$ converge to $\vec 0$ as $y \to + \infty$, we get a fixed point problem for $(U_s, U_p, v_p, V_p)$, which by applying the Contraction Map Theorem is shown to admit a unique solution (in a suitable metric space). 

We then get $\tilde V_s  = V_s  - V_p$, $\tilde v_k = v_k - v_p$ and  $U_k =W - U_p - U_s$: by construction, $U_k$ and $U_s$ satisfy the properties described in the statement of the lemma. 

\noindent
{\sc Step 6} We now provide the details of the argument sketched in {\sc Step 5}. One has to keep in mind that now $W(y)$, $V_s(y)$ and $v_k(y)$ are given functions, and that from~\eqref{e:v} one gets the identities $\tilde V_s  = V_s  - V_p$ and $\tilde v_k = v_k - v_p$. Hence, by integrating the first line in~\eqref{e:s} one gets that the equation for $U_s(y)$
\begin{equation}
\label{e:us}
\begin{split}
          U_s(y) = 
          & \int_{+ \infty}^y   R_{cs} \big(U_s + U^{\sharp}, R_s (V_s -V_p)  (z) , \vec 0 \big) R_s (U_s + U^{\sharp},  V_s - V_p, 0) (z) \Big[ V_s (z) + \\
          & - V_p (z) \Big] dz .\\
\end{split}             
\end{equation}
Combining equations \eqref{e:cs}--\eqref{e:s}, we arrive at the equation for $U_p(y)$
\begin{equation}
\label{e:up}
\begin{split}
          U_p(y) = 
          & \int_{+ \infty}^y  \Big[ R_{cs} \big(W (z), v_k \check r_k (z) + R_s V_s (z) , \vec 0 \big) \check r_k (z) + \\
          & -
          R_{cs}  \big(W (z) -U_s (z) -U_p (z), v_k \check r_{kk} (z), \vec 0 \big) \check r_{kk} (z)
          \Big] v_k (z) dz +\\
          &  
         +  \int_{+ \infty}^y  R_{cs}  \big(W (z) -U_s (z) -U_p (z), v_k \check r_{kk} (z), \vec 0 \big) \check r_{kk} 
          v_p (z) dz + \\
          & + 
          \int_{+ \infty}^y  \Big[ R_{cs} \big(W (z), v_k \check r_k (z) + R_s V_s (z), \vec 0 \big) R_s (z) + \\
          & -
          R_{cs}  \big(U_s (z)+ U^{\sharp}, \tilde V_s R_{ss}, \vec 0 \big) R_{ss} (z)
          \Big] V_s (z) dz +\\
           & +
          \int_{+ \infty}^y  R_{cs}  \big(U_s (z)+ U^{\sharp}, \tilde V_s R_{ss}, \vec 0 \big) R_{ss} (z)
            V_p (z) dz. \\
\end{split}
\end{equation}
To simplify the exposition, we employ the notations $\check r_k= \check r_k (W, v_k, 0)$,
$ \check  r_{kk} = \check r_k (W-U_s -U_p, \tilde v_k, 0)$, $R_s = R_s (W, V_s, 0)$ and $R_{ss} = R_s (U_s + U^{\sharp}, \tilde V_s, 0)$.  
Similarly, we get the equation for $v_p$ 
\begin{equation}
\label{e:vp}
 \begin{split}
         v_p (y) = & \int_{+ \infty}^y     \phi_k \big( W -U_s- U_p, v_k - v_p, \vec 0, 0 )(z) v_p (z) dz \; + \\
&     +    \int_{+\infty}^y \Big[   \phi_k \big( W , v_k, V_s, 0 ) -  \phi_k \big( W -U_s- U_p, v_k - v_p, \vec 0, 0 )(z)
        \Big] v_k (z) dz \;,  \\
  \end{split}      
\end{equation}
and the equation for $V_p$
\begin{equation}
\label{e:Vp}
\begin{split}
         V_p (y) = & \int_0^y e^{\bar \Lambda (y -z)} \Big[ \Lambda (W, v_k, V_s, 0) - 
          \Lambda (U_s +U^{\sharp}, 0, V_s-V_p, 0) \Big] V_s (z) dz \; +  \\
&      +  \int_0^y e^{\bar \Lambda (y -z)} \Big[ \Lambda (W, v_k, V_s, 0) - \bar \Lambda \Big] V_p (z) dz. \\
\end{split}
\end{equation}
Here, we have set
\begin{equation}
\label{e:L}
      \bar \Lambda  = \Lambda \Big(U^{\sharp}, 0, \vec 0, \lambda_k(U^\sharp) \Big), 
\end{equation}      
namely $\bar \Lambda$ is a diagonal matrix whose entries $\lambda_1 (U^{\sharp}), \dots, \lambda_{k-1} (U^{\sharp})$ have all strictly negative real part. 

We now introduce the sets 
\begin{equation*}
\begin{split}   
&           X_s : = \{ U_s \in C^0 ([0, + \infty [ ; \R^n): \; |U_s (y)| \leq C_1 \delta e^{- cy/4}  \; \text{for every} \; y  \} \\
&           X_p^U : = \{ U_p \in C^0 ([0, + \infty [ ; \R^n): \; |U_p (y)| \leq C_2 \delta^2 e^{- cy/4} \; \text{for every} \; y    \} \\ 
&           X_p^V : = \{ V_p \in C^0 ([0, + \infty [ ; \R^{k-1}): \; |V_s (y)| \leq C_3 \delta^2 e^{- cy/4} \; \text{for every} \; y     \} \\ 
&           X_p^v : = \{ v_p \in C^0 ([0, + \infty [ ): \; |v_p (y)| \leq C_4 \delta^2 e^{- cy/4}   \; \text{for every} \; y  \}, \\ 
\end{split}
\end{equation*}
where $c$ is the same as in~\eqref{e:c} and $C_1$, $C_2$, $C_3$ and $C_4$ are suitable constants to be determined. We equip each of the previous spaces with the norm  
$$
    \| F \|_{\ast} : = \sup_{y \in [0, + \infty[} \{ |F(y)| e^{cy/4} \},    
$$
so that $X_s, \; X_p^U, \; X_p^V$ and $X_p^v$ are all complete metric spaces.  We also equip the product space $X_s \times  X_p^U \times  X_p^V \times X_p^v$ with the norm 
$$
    \| (U_s, U_p, V_p, v_p) \|: = \| U_s \|_{\ast} +  \| U_p \|_{\ast} +  \| V_p \|_{\ast} +  \| v_p \|_{\ast}.
$$    
On $X_s \times  X_p^U \times  X_p^V \times X_p^v$ we define the map $T = (T_1, T_2, T_3, T_4)$ by setting $T_1$, $T_2$, $T_3$ and $T_4$ equal to the right hand sides of~\eqref{e:us},~\eqref{e:up},~\eqref{e:Vp} and~\eqref{e:vp} respectively. 

We also recall that from~\eqref{e:uprime} we deduce that $|V_{cs} (y) | \leq \unpo \delta$ for every $y$. Hence, by applying the Local Invertibility Theorem to the map $(V_s, v_k) \mapsto V_{cs}$ defined as in~\eqref{e:vcs}, we obtain the functions $V_s$ and $v_k$ satisfying 
\begin{equation}
\label{e:vsvk}
    |V_s (y)|, \; |v_k(y)| \leq C_5 \delta\;, \qquad \forall y \in [0, + \infty[\;,
\end{equation}
for a suitable constant $C_5 >0$.

We now introduce the following estimate, which is proved in {\sc STEP 7}: there exists a constant $C_6 >0$ such that 
\begin{equation}
\label{e:viesse}
         |V_s (y)| \leq C_6 e^{-c y/2}\delta\;, \qquad \text{for every $y >0$}.   
\end{equation}

By relying on~\eqref{e:viesse} and~\eqref{e:vsvk}, one can then show that, if $\delta$ is small enough, then one can choose the constants $C_1$, $C_2$, $C_3$ and $C_4$ in such a way that $T$ attains values in $X_s \times  X_p^U \times  X_p^V \times X_p^v$ and is a strict contraction. The computations are quite standard and are similar to (but easier than) those performed to prove Lemma 3.9 in~\cite{BiaSpi}, so we skip them. 

Finally, by applying the Contraction Map Theorem we infer that the systems obtained by combining equations~\eqref{e:us},~\eqref{e:up},~\eqref{e:Vp} and~\eqref{e:vp} admit a unique solution belonging to $X_s \times  X_p^U \times  X_p^V \times X_p^v$.

{\sc STEP 7.} We conclude by proving~\eqref{e:viesse}. Let us consider the set  
$$
    X : = \{ V \in C^0 \big( [0, + \infty[, \R^{k-1}\big): \; \| V  \|_s \leq C_6 \delta \},
$$ 
where the value of the constant $C_6$ will be determined in the following and 
$$
    \| V \|_{s} = \sup_{y \in [0, + \infty[} \{  e^{c y /2} |V (y) |\}.
$$
It turns out that $X$, equipped with the distance induced by the norm $\| \cdot \|_{s}$, is a closed metric space. Also, if the constant $C_6$ is big enough, then the map defined by 
$$
    G(V) (y) = e^{\bar \Lambda y} V_s (0) +
    \int_0^y e^{\bar \Lambda (y - z)} 
    \Big[   \Lambda \big(W(z) , v_k(z), V(z), 0 \big) - \bar \Lambda \Big] V (z) dz   
$$
is a strict contraction from $X$ to itself. To prove that $G$ is a contraction, one exploits bounds~\eqref{e:bd} and~\eqref{e:vsvk}. The fixed point of $G$ is the solution of the ODE $V' = \Lambda (W, v_k, V, 0) V$ which attains the value $V_s(0)$ at $y=0$, hence from~\eqref{e:cs} we have $V (y) = V_s(y)$ for every $y$. From the definition of $X$ it follows that $V_s$ satisfies the bound~\eqref{e:viesse}.  
\end{proof}

\begin{remark}
\label{r:in}
As the above proof shows, the heuristic idea underlying Lemma~\ref{l:dec} is that we can ``invert" the construction in~\cite{AnBia} (or in~\cite[pages 45-68]{BiaSpi}). Indeed, in~\cite{AnBia} one first finds $\big(U_k, \tilde{v}_k\big)$ and $\big(U_s, \tilde{V}_s \big)$ by solving~\eqref{e:k} and~\eqref{e:s} respectively and then obtains $(U, v_k, V_s)$ by constructing the perturbation term $(U_p, v_p, V_p)$. Conversely, in Lemma~\ref{l:dec} we are given $(U, v_k, V_s)$ and we want to obtain  $\big(U_k, \tilde{v}_k\big)$, $\big(U_s, \tilde{V}_s \big)$ and $(U_p, v_p,V_p)$. 
\end{remark}

\subsection{Uniqueness of the characteristic wave fan curve}
\label{sus:uwave}
In this subsection, we establish the uniqueness of the characteristic wave fan curve for the boundary characteristic system. It should be noted that the construction of the characteristic wave fan curve is presented in Section~\ref{sus:rie}.

Let $U(t, x) = V(x/t)$ be a self-similar solution of the conservation law~\eqref{e:t:cl} having small enough total variation. Then $V$ has the structure described in Section~\ref{sus:sss} and in particular there exist constant states $U_k^\sharp$, \dots,$U_{n-1}^\sharp\in\R^n$, such that
\be
\label{e:states2}
 V(\xi)=\left\{ \begin{array}{ll}
			U_i^\sharp &    \lambda_i+<\xi<\lambda_{i+1}-,\qquad i=k,\ldots,n-1\\
			U_0 &  \lambda_n+<\xi<+\infty\;,
			\end{array}\right.
\ee
where $\lambda_i+$ and $\lambda_i-$ are defined in~\eqref{e:pm}. Also, in the following we use the same notation as in~\eqref{e:lim}, we denote by $\bar U$ the right limit of $V$ at $0$, namely ${\bar U = \lim_{\xi \to 0^+} V(\xi)}$, and we recall the constants $\delta_1$ and $C$ in the beginning of Section~\ref{s:uni}.
\begin{lemma}
\label{t:u}
         Let the functions $F$ and $B$ satisfy conditions~\eqref{e:eig}--\eqref{e:diss}. Then there exists a sufficiently small constant $\delta$, $ 0<\delta\le\delta_1$, such that, if~\eqref{e:c1} and~\eqref{e:lk1} hold and if $U(t, x) = V(x/t)$ is a self-similar solution of the conservation law~\eqref{e:t:cl} such that 
         \begin{enumerate}
         \item        $TotVar \, V \leq C \delta $,
         \item    all the shocks and the contact discontinuities of $U$ are 
         admissible in the sense of Liu,
         \end{enumerate}
         then the following holds. Let $U_k^{\sharp}$ and $\bar U$ be
         as in~\eqref{e:states2}  and~\eqref{e:lim} respectively and assume that $U_b$ is a state satisfying 
          \begin{enumerate} \setcounter{enumi}{2} 
         \item  there exists a value $\underline U$ such that the following two properties hold:
     \begin{enumerate}
     \item $\underline U$ (on the left) and 
     $\bar U$ (on the right) are connected by a shock (or a contact discontinuity) having speed $0$ and satisfying the Liu admissibility condition;
     \item there exists a solution of~\eqref{e:k}, entirely contained in a ball centered at $(U_k^\sharp,0)$ and having radius $C\delta$. Such a solution satisfies  
     $U_k(0) = U_b$, 
           $\lim_{y \to + \infty} U_k (y) = \underline U$ and $\lim_{y \to + \infty} U'_k (y) = \vec 0$. 
     \end{enumerate}      
                             \end{enumerate}
          Then $U_b =  \tchar (s, U^{\sharp}_k)$ for some $s$, where $\tchar$ is the characteristic wave fan curve defined in~\eqref{e:char1}.   
\end{lemma}
\begin{proof}
The proof exploits many ideas and techniques of Bianchini from~\cite[Theorem $3.2$]{Bia:riemann}. An overview of his analysis can be found in Section~\ref{sus:brie}. The main difference is that we also need to treat the presence of the boundary layers. To deal with this, we employ the monotone convex envelopes of the generalized flux.

The proof is divided into several steps:\\
{\sc Step 1} Fix the orientation of the eigenvector $r_k(U)$ in such a way that 
$$
    \big[  \bar U - U^{\sharp}_k  \big] \cdot r_k (U^{\sharp}_k ) \leq  0, 
$$
where as before the symbol $``\cdot"$ denotes the standard dot product in $\R^n$. By relying on the analysis in \cite{Bia:riemann}, we infer that there exists a value $\bar s<0$ such that the solution of the fixed point problem 
\begin{equation}
\label{e:intermedio}
    \left\{
    \begin{array}{lll}
              \displaystyle{ U(\tau) = 
             U^{\sharp}_k  + \int_{0}^{\tau} \tilde r_k 
              \big( U(z), v (z), \sigma (z)  \big) dz } \\ \\
              \displaystyle{ v(\tau) = 
              f \big(\tau) -  \mathrm{conv}_{[ \bar s, 0 ]} f (\tau) } \\ \\
             \displaystyle{ \sigma (\tau) = \frac{1}{d} \frac{d  \mathrm{conv}_{[\bar s, 0]} f (\tau) }{d \tau}}
    \end{array}
    \right.
\end{equation}
satisfies $U(\bar s) = \bar U$. Here, $ \mathrm{conv}_{[\bar s, 0]}$ denotes the convex envelope of the function $f$ computed on the interval $[\bar s, 0]$ and the function $f$ is defined by
\be 
\label{e:ffeA}
    f(\tau) = \int_{0}^{\tau}  \left[ \phi_k \big( U(z), v(z),\vec 0,  \sigma(z) \big) +d\sigma(z)\right]dz + A,
\ee
where $A$ is an arbitrary constant. Indeed, the solution of~\eqref{e:intermedio} does not depend on the choice of $A$. In~\eqref{e:ffeA}, $\phi_k$ is the same function 
as in~\eqref{e:cs} and the constant $d$ is given by
\be
\label{e:kappa}
    d : = - \Bigg[ \, \frac{\partial \phi_k (U, v_k, V_s, \sigma) }{ \partial \sigma} \Big|_{U = U^{\sharp}_k, v_k =0, V_s = \vec 0, \sigma= \lambda_k (U^{\sharp}_k) } \Bigg] 
\ee
and it is strictly positive by construction (see~\cite[formulas (3.32)-(3.33)]{BiaSpi}). 
 \\
{\sc Step 2} By combining the assumption (3a) in the statement of the lemma with the analysis in~\cite[Section 3]{Bia:riemann}, we deduce that there are two cases, the first one being that there exists $\underline s < \bar s$ such that the fixed point problem 
\begin{equation}
\label{e:0speed}
    \left\{
    \begin{array}{lll}
              \displaystyle{ U(\tau) = 
             \bar U + \int_{\bar s}^{\tau} \tilde r_k 
              \big( U(z), v (z), \sigma (z)  \big) dz } \\ \\
              \displaystyle{ v(\tau) = 
              f (\tau) -  \mathrm{conv}_{[ \underline s, \bar s ]} f (\tau) } \\ \\
             \displaystyle{ \sigma (\tau) = \frac{1}{d} \frac{d  \mathrm{conv}_{[\underline s, \bar s]} f (\tau) }{d \tau}}
    \end{array}
    \right.
\end{equation} 
defined on the interval $[\underline s, \bar s]$ with
\be 
\label{e:ffe}
    f(\tau) = \int_{\bar s}^{\tau}   \left[ \phi_k \big( U(z), v(z),\vec 0,  \sigma(z) \big) +d\sigma(z)\right]dz
\ee
has a solution $(U,v,\sigma)$ that satisfies $U(\underline s) = \underline U$.

The second case is that there exists a $\underline s > \bar s$ such that the solution of the fixed point problem 
\begin{equation}
\label{e:wrong}
    \left\{
    \begin{array}{lll}
              \displaystyle{ U(\tau) = 
             \bar U + \int_{\bar s}^{\tau} \tilde r_k 
              \big( U(z), v (z), \sigma (z)  \big) dz } \\ \\
              \displaystyle{ v(\tau) = 
              \ell(\tau) -  \mathrm{conc}_{[  \bar s, \underline s ]} \ell (\tau) } \\ \\
             \displaystyle{ \sigma (\tau) = \frac{1}{d} \frac{d  \mathrm{conc}_{[ \bar s, \underline s]} \ell (\tau) }{d \tau}}
    \end{array}
    \right.
\end{equation} 
satisfies $U(\underline s) = \underline U$. Note that the fixed point problem is now defined on the interval $[ \bar s, \underline s ]$. In the previous expression, the function $\ell$ is defined as in~\eqref{e:ffe}, namely 
\be 
\label{e:ell}
          \ell (\tau) = \int_{\bar s}^{\tau}  \left[ \phi_k \big( U(z), v(z),\vec 0,  \sigma(z) \big) +d\sigma(z)\right]dz.
\ee
Also, $\mathrm{conc}_{[ \bar s, \underline s]} \ell$ denotes the concave envelope of $\ell$ on $[ \bar s, \underline s]$, namely $\mathrm{conc}_{[ \bar s, \underline s]} \ell : = - \mathrm{conv}_{[ \bar s, \underline s]} (-\ell)$. In {\sc Step 3}, this second possibility is ruled out. \\
{\sc Step 3} The goal in this step is to show that the second case presented above cannot occur. To prove this, we proceed by contradiction and exploit Lemma $3.1$ in~\cite{Bia:riemann}. 

Assume that~\eqref{e:wrong} defines a zero speed shock (or contact discontinuity) which is admissible in the sense of Liu and which connects $\bar U$ on the right to $\underline U$ on the left. In other words, assume that the solution $(U,v,\sigma)$ of~\eqref{e:wrong} satisfies $U(\underline s) = \underline U$, $\mathrm{conc}_{[ \bar s, \underline s]} \ell  \equiv 0$ and 
\be
\label{e:ieffe} 
         \ell (\tau) \leq 0 \quad \text{ for every $\tau \in [\bar s, \underline s]$}.
\ee
Indeed, the fact that $\mathrm{conc}_{[ \bar s, \underline s]} \ell (\tau)\equiv 0$ is equivalent to the admissibility of the shock. 
\begin{figure}
\begin{center}
\caption{The case we rule out in {\sc Step 3}} 
\psfrag{a}{$\bar s$}  \psfrag{f}{$f (\tau)$} \psfrag{b}{$\bar s + h$}
\psfrag{c}{$\underline s$}   \psfrag{g}{$\ell^h(\tau)$} \psfrag{h}{$\ell$} 
\label{f:e}
\includegraphics[scale=0.6]{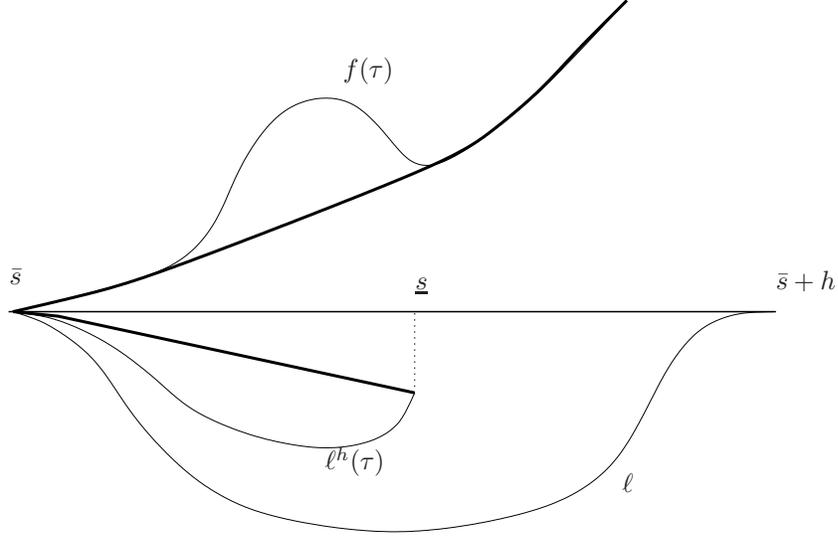}
\end{center}
\end{figure}
Let us now fix $h< \min \{ | \bar s  |, |\underline s - \bar s | \}$ (see Figure~\ref{f:e}) and consider the map 
\begin{equation}
\label{e:acca}
    \left\{
    \begin{array}{lll}
              \displaystyle{ T_1 [U, v, \sigma](\tau)= 
             \bar U + \int_{\bar s}^{\tau} \tilde r_k 
              \big( U(z), v (z), \sigma (z)  \big) dz } \\ \\
              \displaystyle{ T_2 [U, v, \sigma](\tau) = 
              \ell (\tau)-  
              \mathrm{conc}_{[  \bar s, \bar s + h ]}  \ell (\tau) } \\ \\
             \displaystyle{ T_1 [U, v, \sigma](\tau) = \frac{1}{d} \frac{d  \mathrm{conc}_{[ \bar s, \bar s + h]}  \ell (\tau) }{d \tau}},
    \end{array}
    \right.
\end{equation} 
defined on the interval $[\bar s, \bar s + h]$. Also, let us define the norm 
\be 
\label{e:hnorm}
    \| (U, v, \sigma) \|_h: = \zeta \| U \|_{C^0 ([ \bar s, \bar s + h]) } +
     \zeta \| v\|_{C^0( [ \bar s, \bar s + h] )} + \| \sigma \|_{L^1 ([ \bar s, \bar s + h])},
\ee
where $\zeta$ is a small enough positive constant to be determined.
Then, we get that, if $\zeta$ is small enough, the map $T= (T_1, T_2, T_3)$ defined in~\eqref{e:acca} is a strict contraction from a suitable metric space into itself with respect to the norm $\| \cdot \|_h$. In the following, we denote by $\gamma^h = (U^h, v^h, \sigma^h)$ the unique fixed point of $T$ and by $ \ell^h$ the quantity 
$$
      \ell^h (\tau) =  \ell (U^h, v^h, \sigma^h) (\tau)\; .
$$
By relying on the same argument as in the proof of Lemma 3.1 in~\cite{Bia:riemann}, we get that 
\be 
\label{e:isigma}
         \sigma^h (\tau) \leq 0 \quad \text{for every $\tau \in [ \bar s, \bar s + h]$}. 
\ee
Since~\eqref{e:ell} implies that $l(\bar s)=0$, then $\mathrm{conc}_{[ \bar s, \bar s + h]}  \ell^h (\bar s)=0$. Hence inequality~\eqref{e:isigma} implies that 
\be
\label{e:ieffeacca}
    \ell^h (\tau) \leq  \mathrm{conc}_{[ \bar s, \bar s + h]}  \ell^h (\tau) \leq 0 
   \quad \text{for every $\tau \in [\bar s, \bar s + h]$}.
\ee
Let us now return to~\eqref{e:intermedio}: since $\sigma (\tau) >0$ for $\tau\in]\underline s,0]$, then by arguing as before we get $$0 <  \mathrm{conv}_{[ \bar s, 0]} f (\tau)\leq f (\tau)$$ for every $\tau \in ]\bar s, 0]$. Hence~\eqref{e:ieffeacca} implies
\be
\label{e:contra}
         \ell^h (\tau ) < f (\tau) \quad \text{for every $\tau \in \,  ]\bar s, \bar s + h]$}.
\ee
Here and in the following, we denote by $f$ the restriction to the interval $[\bar s, \bar s + h]$ of the same function $f$ as in~\eqref{e:intermedio}. Next, we estimate
$\| f - \ell^h \|_{C^0([ \bar s, \bar s + h])}$. First, we point out that, for every $\tau \in [ \bar s, \bar s + h]$, we have
\begin{equation}
\label{e:effeacca}
\begin{split}
          f (\tau) & - \ell^h (\tau) 
       = \int_{\bar s}^{\tau} \Big[ \phi_k (U, v, \sigma)  - \phi_k( U^h, v^h,   \sigma^h) \Big] (\tau) d \tau + \int_{\bar s}^{\tau} d(\sigma-\sigma^h)d\tau\\
&       \leq \unpo h \| U - U^h \|_{C^0([ \bar s, \bar s + h])} + \unpo h
          \| v - v^h \|_{C^0[ \bar s, \bar s + h]} +
          \zeta \| \sigma - \sigma^h \|_{L^1([ \bar s, \bar s + h])}.
\end{split}
\end{equation}
Here and in the sequel, we denote by $\unpo$ any constant depending only on the flux $F$ and on the state $U^{\sharp}_k$. Also, we exploit the estimate ${|\partial \phi_k / \partial \sigma +d| \leq \unpo \delta}$ and, by taking $\delta$ sufficiently small, we can assume ${|\partial \phi_k / \partial \sigma +d| \leq \zeta^2}$. By choosing $h\le\zeta^3$ in~\eqref{e:effeacca}, relying on the definition~\eqref{e:hnorm} and using the fact that the map $T = (T_1, T_2, T_3)$ defined in~\eqref{e:acca} is Lipschitz continuous with constant $1/2$ and that $\gamma^h$ is the fixed point of $T$, we eventually get 
\be 
\label{e:econtra}
      f (\tau) - \ell^h (\tau)  \leq \unpo \zeta^2 \| \gamma - \gamma^h  \|_h \leq \frac{\zeta^2 }{1 - 1/2}
      \| \gamma - T (\gamma) \|_h, 
\ee
where $\gamma = (U, v, \sigma)$ is the restriction of the solution of~\eqref{e:intermedio} to the interval $[ \bar s, \bar s + h]$. To estimate $\| f - f^h \|_{C^0([ \bar s, \bar s + h])}$, it remains to provide a bound on $ \| \gamma - T (\gamma) \|_h$. By recalling that the solution of~\eqref{e:intermedio} satisfies $U(\bar s ) = \bar U$ and taking into account~\eqref{e:acca}, we deduce 
\be
\label{e:g1}
         U  - T_1 (\gamma) \equiv 0 \quad \text{on $[\bar s, \bar s + h]$}
\ee
and 
$$
     | v(\tau) - T_2 (\gamma) (\tau) | = 
     | \mathrm{conc}_{[  \bar s, \bar s + h ]} f (\tau) -
     \mathrm{conv}_{[  \bar s, 0 ]} f (\tau) |. 
$$
By relying on the chain of inequalities 
\be 
\label{e:g2a}
    \max_{\tau \in [\bar s , \bar s+ h] } f \ge \mathrm{conc}_{[  \bar s, \bar s + h ]} f (\tau) \ge f (\tau) \ge \mathrm{conv}_{[  \bar s, 0 ]} f (\tau) > 0 \ge \ell^h (\tau) \quad  \forall \, \tau \in ]\bar s , \bar s+ h] ,
\ee
we obtain
\be
\label{e:g2}
\begin{split}
           | v(\tau) - T_2 (\gamma) (\tau) | & \leq \mathrm{conc}_{[  \bar s, \bar s + h ]} f (\tau) -
     \mathrm{conv}_{[  \bar s, 0 ]} f (\tau)  \leq  \max_{\tau \in [\bar s , \bar s+ h] } f  \\
&     \leq \| f - \ell^h \|_{C^0( [\bar s , \bar s+ h] )}. \\
\end{split}
\ee
Finally, let us consider the quantity
\be 
\label{e:g3a}
    \| \sigma - T_3 (\gamma) \|_{L^1 ([\bar s , \bar s+ h])} = \int_{\bar s}^{  \bar s+ h }
    \left| \frac{  \mathrm{conv}_{[  \bar s, 0 ]} f  }{d \tau }   - \frac{\mathrm{conc}_{[  \bar s, \bar s + h ]} f }{d \tau }\right| d \tau\;.
\ee
Since $ \mathrm{conv}_{[  \bar s, 0 ]} f / d \tau - \mathrm{conc}_{[  \bar s, \bar s + h ]} f / d \tau$ is non-decreasing, the set of zeroes is an interval, which can be empty or reduced to a point. Let us fix $s^{\ast}$ such that 
$$
    \mathrm{conv}_{[  \bar s, 0 ]} f / d \tau - \mathrm{conc}_{[  \bar s, \bar s + h ]} f / d \tau < 0 
     \quad \text{in $[\bar s, s^{\ast}[$}
$$ 
and 
$$
   {\mathrm{conv}_{[  \bar s, 0 ]} f / d \tau - \mathrm{conc}_{[  \bar s, \bar s + h ]} f / d \tau \ge 0}
   \quad \text{in $[s^{\ast}, \bar s + h]$}.
$$     
Using~\eqref{e:g3a}, we estimate
\be 
\label{e:g3}
\begin{split}
            \| \sigma - T_3 (\gamma) \|_{L^1 ([\bar s , \bar s+ h])}
 & =        
       \int_{\bar s}^{   s^{\ast} } \Big[
      \mathrm{conc}_{[  \bar s, \bar s + h ]} f / d \tau - \mathrm{conv}_{[  \bar s, 0 ]} f / d \tau \Big] d \tau +  \\
& \quad + \int_{s^{\ast}}^{\bar s + h}
      \Big[ \mathrm{conv}_{[  \bar s, 0 ]} f / d \tau - \mathrm{conc}_{[  \bar s, \bar s + h ]} f / d \tau \Big] d \tau  = \\
&   = 2 \Big[  \mathrm{conc}_{[  \bar s, \bar s + h ]} f  - 
         \mathrm{conv}_{[  \bar s, 0 ]} f  \Big]  (s^{\ast})  + \\
&  \quad +   \Big[ \mathrm{conv}_{[  \bar s, 0 ]} f - \mathrm{conc}_{[  \bar s, \bar s + h ]} f \Big]  (\bar s + h)   \leq \\
&      \leq 2 \Big[  \mathrm{conc}_{[  \bar s, \bar s + h ]} f  - 
         \mathrm{conv}_{[  \bar s, 0 ]} f  \Big]  (s^{\ast}) \leq   \\
&      \leq 2 \| f - \ell^h \|_{C^0( [\bar s , \bar s+ h] )}.          \\
\end{split}
\ee
To get the last two inequalities, we exploit~\eqref{e:g2a} and~\eqref{e:g2} respectively. 

Combining~\eqref{e:g1},~\eqref{e:g2} and~\eqref{e:g3} with~\eqref{e:econtra}, we arrive at 
$$
    \| f (\tau) - \ell^h (\tau) \|_{C^0 ([\bar s, \bar s + h])} \leq \unpo \zeta 
      \| f (\tau) - \ell^h (\tau) \|_{C^0 ([\bar s, \bar s + h])}   .
$$
By taking $\zeta$ sufficiently small, the above inequality implies $ \| f (\tau) - \ell^h (\tau) \|_{C^0 ([\bar s, \bar s + h])} =0$. Since this contradicts~\eqref{e:contra}, we conclude that~\eqref{e:wrong} cannot hold. \\
{\sc Step 4} From here and on, we investigate the first case stated in {\sc Step 2}, namely the fixed point problem~\eqref{e:0speed}--\eqref{e:ffe}. In this step, we ``glue together"~\eqref{e:intermedio} and~\eqref{e:0speed} by finding a fixed point problem defined on the interval $[\underline s, \bar s]$. 

First, we point out that the function $f(\tau)$ defined by~\eqref{e:ffeA} for $\tau\in[\bar s,0]$ and by \eqref{e:ffe} for $\tau\in[\underline s,\bar s[$ is continuous on $[\underline s,0]$, because in~\eqref{e:ffeA} the constant $A$ is chosen in such a way that $f(\bar s)=0$. We define the function $g$ on $[\underline s, 0]$ by setting
\be
\label{e:g}
         g(\tau) =
         \left\{
         \begin{array}{ll}
                    \mathrm{conv}_{[  \bar s, 0]} f (\tau)  & \text{if $\tau \in [ \bar s, 0]$} \\ \\
            \mathrm{conv}_{[ \underline s, \bar s]} f (\tau)  &
             \text{if $\tau \in [\underline s, \bar s [$} .\\
         \end{array}
         \right.
\ee
Since we assume that $\underline U$ and $\bar U$ are connected by a shock (or a contact discontinuity) with speed zero, then $d  \,  \mathrm{conv}_{[ \underline s, \bar s ]} f / d \tau$ is actually $0$ on $[\underline s, \bar s ]$. Also, by assumption $d \, \mathrm{conv}_{[  \bar s, 0]} f / d \tau$ is strictly positive on $[\bar s, 0]$. It then turns out that the function $g$ is convex. In addition, by construction $g (\tau) \leq  f (\tau)$ in $[\underline s, 0]$, hence
\be
\label{e:d1}
        g(\tau) \leq  \mathrm{conv}_{[  \underline s, 0]} f (\tau) \; \text{for every $\tau \in [\underline s, 0]$}.
\ee
Conversely, the restriction of $\mathrm{conv}_{[  \underline s, 0]} f (\tau)$ to the interval $[\bar s, 0]$ is a convex function less than or equal to $f$, hence 
$$
     \mathrm{conv}_{[  \underline s, 0]} f (\tau) \leq  \mathrm{conv}_{[  \bar s, 0]} f (\tau) \quad  \text{for every $\tau \in [\bar s, 0]$}.
$$
By applying the same argument to the interval $[\underline s, \bar s ]$, we deduce that 
\be
\label{e:d2}
        g(\tau) \ge  \mathrm{conv}_{[  \underline s, 0]} f (\tau) \quad \text{for every $\tau \in [\underline s, 0]$}
\ee
and therefore, by combining inequality~\eqref{e:d1} with~\eqref{e:d2}, we arrive at the identity
$$
    g(\tau) =  \mathrm{conv}_{[  \underline s, 0]} f (\tau) \quad 
     \text{for every $\tau \in [\underline s, 0]$}.
$$
We finally deduce that the solution of the fixed point problem (defined on the interval $[\underline s, 0]$) 
\begin{equation}
\label{e:0+}
    \left\{
    \begin{array}{lll}
              \displaystyle{ U(\tau) = 
             U^{\sharp}_k + \int_{\underline s}^{\tau} \tilde r_k 
              \big( U(z), v (z), \sigma (z)  \big) dz } \\ \\
              \displaystyle{ v(\tau) = 
              f (\tau)-  \mathrm{conv}_{[ \underline s, 0]} f (\tau) } \\ \\
             \displaystyle{ \sigma (\tau) = \frac{1}{d} \frac{d  \mathrm{conv}_{[\underline s, 0]} f (\tau) }{d \tau}}
    \end{array}
    \right.
\end{equation} 
satisfies $U(\bar s)  = \bar U$ and $U(\underline s) = \underline U$. \\
{\sc Step 5}
We now show that, under the hypotheses of the lemma, one can use the monotone convex 
envelope in~\eqref{e:0+}  instead of the convex one. The monotone convex envelope is defined in~\eqref{e:monconv}. Let $C^{1, 1}_L([\underline s, 0]) $ denote the space of the functions $f$ defined on $[\underline s, 0]$ such that $f$ is continuously differentiable functions and $f'$ is Lipschitz continuous with Lipschitz constant bounded by $L$. If $f \in C^{1, 1}_L ([\underline s, 0]) $, then $\mathrm{conv}_{[\underline s, 0]} f \in C^{1, 1}_L ([\underline s, 0]) $. This follows from more general results discussed  by Griewank and Rabier in~\cite{GriewankRabier}. Also, if $f  \in C^{1, 1}_L  ([\underline s, 0]) $, then 
\begin{equation}
\label{e:mon}
     \mathrm{monconv} f(\underline s, 0]) = 
     \left\{
     \begin{array}{ll}
                  \mathrm{conv} f ([\underline s, 0]) (\tau) & \text{if $\tau \ge \tau_0$}       \\
                   \mathrm{conv} f ([\underline s, 0]) (\tau_0) & \text{if $\tau \leq \tau_0$}, \\    
     \end{array}
     \right.
\ee
where 
$$
    \tau_0 : = \min \left\{ \tau \in [\underline s, 0]: \;  \frac{d}{d \tau}  \mathrm{conv} f ([\underline s, 0]) \ge 0 \right\}.
$$
If the derivative of the function $ \mathrm{conv} f ([\underline s, 0])$ is strictly negative in $[\underline s, 0]$, then one sets $\tau_0 = \underline s$ and then, relation~\eqref{e:mon} is still valid. The proof of~\eqref{e:mon} can be found, for example, in~\cite[Page 36]{BiaSpi}.
Since in~\eqref{e:0+} $\sigma = \dfrac{1}{d}\dfrac{d}{d\tau}\mathrm{conv} f ([\underline s, 0]) \ge 0 $, then by exploiting~\eqref{e:mon}, we infer that instead of the fixed point problem~\eqref{e:0+}, we can consider the following one:
\begin{equation}
\label{e:fpmoncon}
    \left\{
    \begin{array}{lll}
              \displaystyle{ U(\tau) = 
             U^{\sharp}_k  + \int_{0}^{\tau} \tilde r_k 
              \big( U(z), v (z), \sigma (z)  \big) dz } \\ \\
              \displaystyle{ v(\tau) = 
              f (\tau) -  \mathrm{monconv}_{[ \underline s, 0 ]} f (\tau) } \\ \\
             \displaystyle{ \sigma (\tau) = \frac{1}{d} \frac{d  \mathrm{monconv}_{[\underline s, 0]} f (\tau) }{d \tau}}\;,
    \end{array}
    \right.
\end{equation}
which is still defined on the interval $[\underline s, 0]$.\\
{\sc Step 6} We now exploit assumption (3b) in the statement of the lemma. Namely, there exists a solution $(U, v_k)$ of the system  
\be
\label{e:uv}
     \left\{
         \begin{array}{llll}
                  U'= R_{cs} (U, v_k \check{r}_k, 0) 
                  \check{r}_k (U, v_k, \vec 0,  0)  v_k  \\
                  v_k'= \phi_k ( U,  v_k, \vec 0, 0) v_k \\
         \end{array}
         \right.
\ee
which satisfies $U(0) = U_b$, $\lim_{y \to + \infty} U(y) = \underline U$ and $\lim_{y \to + \infty} U'(y) = \vec 0$. Recall that $R_{cs}$ and $\phi_k$ are constructed in the proof of Lemma~\ref{l:dec}.

If $\underline U = U_b$, then the conclusion of the lemma follows from the previous steps. So in the following we assume $\underline U \neq U_b$.  Since any point of the form $(U, 0)$ is an equilibrium for system~\eqref{e:uv}, then $v_k (y) \neq 0$ for every $y$. 

Let us first consider the case when $v_k (y) >0$ for every $y$. The change of variables 
$$
    \frac{d \tau}{dy } = v_k (y)
$$
maps the interval $[0, + \infty[$ into some interval $[s, \underline s[$ and $(U, v_k)$ satisfies
$$
    \left\{
         \begin{array}{llll}
                  \displaystyle{U (\tau) = \underline U+ \int_{\underline s}^{\tau} \tilde r_k (U,  v_k, \vec 0,  0)  dz }  \\
                 \displaystyle{ v_k (\tau) = \int_{\underline s}^{\tau} \phi_k ( U,  v_k, \vec 0, 0)  dz. } \\
         \end{array}
         \right.
$$
In the previous expression, we employ the equalities $\tilde r_k = R_{cs} \check{r}_k$ and ${\lim_{y \to + \infty} v_k (y) =0}$. The last condition follows from the fact that $U'(y)$ converges to $0$ and $|\tilde r_k|$ is bounded away from $0$. By setting 
$$
    f(\tau)=  \int_{\underline s}^{\tau} \phi_k ( U(z),  v_k(z), \vec 0, 0)  dz  ,
$$
we have $v_k(\tau) = f(\tau) \ge 0$ for every $\tau\in[s,\underline s[$. Since $f(0)=0$, then $0 \equiv \mathrm{monconv}_{[s, \underline s[} f$ and hence we can rewrite the previous system as 
\be
\label{e:fpbl}
\left\{
    \begin{array}{lll}
              \displaystyle{ U(\tau) = 
              \underline U  + \int_{\underline s}^{\tau} \tilde r_k 
              \big( U(z), v (z), \sigma (z)  \big) dz } \\ \\
              \displaystyle{ v(\tau) = 
              f (\tau)-  
              \mathrm{monconv}_{[s, \underline s[}  f (\tau) } \\ \\
             \displaystyle{ \sigma (\tau) = \frac{1}{d} \frac{d  \mathrm{monconv}_{[s, 
             \underline  s[   }  f (\tau) }{d \tau}=0\;.} \\
    \end{array}
    \right.
\ee
By arguing as in {\sc Step 4}, we can then ``glue together"~\eqref{e:fpmoncon} and~\eqref{e:fpbl} and obtain that the solution of the fixed point problem 
$$
    \left\{
    \begin{array}{lll}
              \displaystyle{ U(\tau) = 
             U^{\sharp}_k  + \int_{0}^{\tau} \tilde r_k 
              \big( U(z), v (z), \sigma (z)  \big) dz } \\ \\
              \displaystyle{ v(\tau) = 
              f (\tau) -  \mathrm{monconv}_{[  s, 0 ]} f (\tau) } \\ \\
             \displaystyle{ \sigma (\tau) = \frac{1}{d} \frac{d  \mathrm{monconv}_{[ s, 0]} f (\tau) }{d \tau}},
    \end{array}
    \right.
$$
for $\tau\in[s,0]$ satisfies $U(s) = U_b$. Following a similar argument as in {\sc Step 3}, we can rule out the possibility that $v_k (y) < 0$ in~\eqref{e:uv}.  The proof of the lemma is complete.
\end{proof}

\subsection{Conclusion of the proof}\label{sus:conchar}
Here, we complete the proof of the theorem. In Subsection~\ref{sss:char}, the boundary characteristic case is treated and in Subsection~\ref{sss:nonchar}, the non-characteristic.

\subsubsection{Boundary characteristic case}
\label{sss:char}
By relying on Lemma~\ref{l:sfs}, we infer that the self-similar function $V$ satisfies
\be
\label{e:states3}
V(\xi)=\left\{ \begin{array}{ll}
			U_i^\sharp &    \lambda_i+<\xi<\lambda_{i+1}-,\qquad i=k,\ldots,n -1 \\
			U_n^\sharp=U_0 &  \lambda_n+<\xi<+\infty\;.
			\end{array}\right.
\ee
By relying on assumptions (3) and (4) in the statement of Theorem~\ref{t:ui}, we have that, for any $i=k,\dots,n-1$, $U_i^\sharp$ and $U_{i+1}^\sharp$ are connected by a sequence of rarefactions and shocks (or contact discontinuities) that are admissible in the sense of Liu. Hence, by exploiting Theorem 3.1 in~\cite{Bia:riemann}, we get that, for every $i=k,\dots, n-1 $,  
$$
    U_{i}^{\sharp} = T_{i} (s_{i}, U^{\sharp}_{i+1 })
$$  
for some $s_i$ small. In the previous expression, $T_i$ denotes the same $i$-wave fan curve as in formula~\eqref{e:ciao}. Hence,  
\be 
\label{e:c:bu}
    U^{\sharp}_k = T_{k+1} \Big( s_{k+1}, T_{k+2} (\dots, T_n( s_n, U_0) \dots ) \Big)
\ee
for a suitable vector $(s_{k+1}, \dots s_n)$. 

We now exploit assumption (5) in the statement of Theorem~\ref{t:ui} and we apply Lemma~\ref{l:dec} to the function $W$. Let $U_k(y)$ be the same function as in the statement of Lemma~\ref{l:dec}, then by setting $U_k (0) = U_b$ and applying Lemma~\ref{t:u} we obtain $U_b = \tchar (U^{\sharp}_k, s_k)$ for some small $s_k$. By comparing the construction in~\cite[Section 3.2.3]{BiaSpi}
with the proof of Lemma~\ref{l:dec}, one infers that $U_D = \psi (s_1, \dots, s_k, U^{\sharp}_k)$, for a suitable vector $(s_1, \dots ,s_k)$, where $\psi$ is the same function as in~\eqref{e:psi}. Hence, relation~\eqref{e:c:bu} leads to 
$$
     U_D = \psi \Big( s_1, \dots, s_k,    T_{k+1} 
         \big( s_{k+1}, T_{k+1} (\dots, T_n( s_n, U_0) \dots ) \big)  \Big)
         = \phi(s_1, \dots ,s_n, U_0).
$$
The analysis in~\cite[page 67-68]{BiaSpi} ensures that the map $\phi$ is locally invertible, hence given $U_D$ and $U_0$ such that $|U_D -U_0|$ is small enough, the value of $(s_1, \dots ,s_n)$ is uniquely determined by the above relation. 

In conclusion, if $U$ is a solution satisfying all the assumptions of Theorem~\ref{t:ui}, then the value $U(t, x)$ is uniquely determined for a.e. ${(t, x) \in [0, + \infty[ \times [0, + \infty[}$ and it can be obtained as follows. Set as in~\eqref{e:states3} $U^{\sharp}_n=U_0$ and for any $i=k, \dots n-1$ define inductively $U^{\sharp}_i $ as $T_{i+1}(s_{i+1}, U^{\sharp}_{i+1})$. Also, let $(U_i, v_i, \sigma_i)$ be the fixed point of~\eqref{e:biabre} (if $s_i<0$, otherwise one takes the concave envelope instead of the convex one). Also, assume for now that $s_k<0$ and denote by $(U_k,v_k,\sigma_k)$ the fixed point of \eqref{e:char1} and by $\bar s$ the value
$$
\bar s=\max\{\tau\in[s_k,0] \;:\;\sigma_k(\tau)=0\}.
$$
Then,
$$
    U(t, x) = 
    \left\{
    \begin{array}{ll}             
    		U_k(\tau) & \text{if $x/t= \sigma_k (\tau)$\quad $\tau\in[\bar s,0]$}\\
               U^{\sharp}_i  &  \text{if $ \sigma_i (0) < x / t < \sigma_{i+1}(s_{i+1}) $}\quad\text{for}\; i=k,\dots,n-1\\     
               U_j (\tau)       &   \text{if  $x/t =  \sigma_j (\tau)$ } \quad\text{for}\; j=k+1,\dots,n\\
               U^{\sharp}_n &    \text{if  $x/t > \sigma_n (0).$}
    \end{array}
    \right.
$$
If $s_k>0$, one needs to take in~\eqref{e:char1}, the monotone concave envelope and define the value $\bar s$ as
$$
\bar s=\min\{\tau\in[0,s_k] \;:\;\sigma_k(\tau)=0\}\;.
$$
The proof for the boundary characteristic case is thus complete.

\subsubsection{Non characteristic boundary case}\label{sss:nonchar}
We now focus on the case when the boundary is non characteristic, namely we assume that~\eqref{e:nchar} holds.

By arguing as in Section~\ref{sss:char} and relying on assumptions (1)-(4) in the statement of the theorem, we infer that $\bar U$, the trace of $U$ on the $t$ axis, satisfies 
\be 
\label{e:nc:bu}
    \bar U = T_{n-p+1} \Big( s_{n-p+1}, T_{n-p+2} (\dots, T_n( s_n, U_0) \dots ) \Big)
\ee
for a suitable vector $(s_{n-p+1}, \dots s_n)$. 

Let us now exploit the assumption (5) in the statement of the theorem: since the boundary is non characteristic, then $f$ is locally invertible and from the relations $f(\bar U) = f(\underline U) $ and $|\bar U - \underline U| \leq C \delta$, we obtain $\bar U = \underline U$. Hence, system~\eqref{e:i:bl} becomes
$$
\left\{
\begin{array}{ll}
          B(W) W' = F(W) - F(\bar U) \\
          W(0) = U_D \qquad \lim_{y \to + \infty } W (y) = \bar U
\end{array}
\right.
$$
and by exploiting again the fact that $DF (\bar U)$ is invertible, we get that $U_D$ lies on the stable manifold of the above system about the equilibrium point $\bar U$, namely
\be
\label{e:nc:ud}
   U_D = \varphi \Big( s_1, \dots, s_{n-p}, \bar U \Big)
\ee
for a suitable vector $(s_1, \dots, s_{n-p})$ and for the same function $\varphi$ as in~\eqref{e:ud}. 

By combining~\eqref{e:nc:bu} with~\eqref{e:nc:ud}, we arrive at
$$
    U_D = \varphi \Bigg( s_1, \dots, s_{n-p}, T_{n-p+1} \Big( s_{n-p+1}, T_{n-p+1} (\dots, T_n( s_n, U_0) \dots ) \Big) \Bigg).
$$
Once $U_D$ and $U_0$ are given, this relation uniquely determines $(s_1, \dots ,s_n)$ and, thus, the function $U$ satisfying assumptions (1)--(5) in the statement of the theorem. \qed

 \section*{Acknowledgements}       
The authors would like to thank Professor Constantine Dafermos for proposing this project and Professors Stefano Bianchini and Athanasios Tzavaras for valuable discussions.
Christoforou was partially supported by the National Sciences Foundation under the grant 
DMS 0803463 and the Texas Advanced Research Program under the grant 003652--0010--2007. 
\bibliography{biblio_t}

\def\cprime{$'$}
\begin{thebibliography}{10}

\bibitem{AmbrosioFuscoPallara}
L.~Ambrosio, N.~Fusco, and D.~Pallara.
\newblock {\em Functions of bounded variation and free discontinuity problems}.
\newblock Oxford Mathematical Monographs. The Clarendon Press Oxford University
  Press, New York, 2000.

\bibitem{AnBia}
F.~Ancona and S.~Bianchini.
\newblock Vanishing viscosity solutions of hyperbolic systems of conservation
  laws with boundary.
\newblock In {\em ``{WASCOM} 2005''---13th {C}onference on {W}aves and
  {S}tability in {C}ontinuous {M}edia}, pages 13--21. World Sci. Publ.,
  Hackensack, NJ, 2006.

\bibitem{Bia:riemann}
S.~Bianchini.
\newblock On the {R}iemann problem for non-conservative hyperbolic systems.
\newblock {\em Arch. Ration. Mech. Anal.}, 166(1):1--26, 2003.

\bibitem{BiaBrevv}
S.~Bianchini and A.~Bressan.
\newblock Vanishing viscosity solutions of nonlinear hyperbolic systems.
\newblock {\em Ann. of Math. (2)}, 161(1):223--342, 2005.

\bibitem{BiaSpi}
S.~Bianchini and L.V. Spinolo.
\newblock The boundary {R}iemann solver coming from the real vanishing
  viscosity approximation.
\newblock {\em Arch. Ration. Mech. Anal.}, 191(1):1--96, 2009.

\bibitem{Bre:book}
A.~Bressan.
\newblock {\em {Hyperbolic systems of conservation laws. The one-dimensional
  Cauchy problem}}, volume~20 of {\em Oxford Lecture Series in Mathematics and
  its Applications}.
\newblock Oxford University Press, Oxford, 2000.

\bibitem{Bre:notes}
A.~Bressan, D.~Serre, M.~Williams, and K.~Zumbrun.
\newblock {\em Hyperbolic systems of balance laws}, volume 1911 of {\em Lecture
  Notes in Mathematics}.
\newblock Springer, Berlin, 2007.
\newblock Lectures given at the C.I.M.E. Summer School held in Cetraro, July
  14--21, 2003. Edited and with a preface by Pierangelo Marcati.

\bibitem{D1}
C.~M. Dafermos.
\newblock Solution of the {R}iemann problem for a class of hyperbolic systems
  of conservation laws by the viscosity method.
\newblock {\em Arch. Rational Mech. Anal.}, 52:1--9, 1973.

\bibitem{D}
C.~M. Dafermos.
\newblock {\em Hyperbolic conservation laws in continuum physics}, volume 325
  of {\em Grundlehren der Mathematischen Wissenschaften [Fundamental Principles
  of Mathematical Sciences]}.
\newblock Springer-Verlag, Berlin, second edition, 2005.

\bibitem{DubLeF}
F.~Dubois and P.~LeFloch.
\newblock Boundary conditions for nonlinear hyperbolic systems of conservation
  laws.
\newblock {\em J. Differential Equations}, 71(1):93--122, 1988.

\bibitem{Gis}
M.~Gisclon.
\newblock \'{E}tude des conditions aux limites pour un syst\`eme strictement
  hyperbolique, via l'approximation parabolique.
\newblock {\em J. Math. Pures Appl. (9)}, 75(5):485--508, 1996.

\bibitem{GisSerre}
M.~Gisclon and D.~Serre.
\newblock \'{E}tude des conditions aux limites pour un syst\`eme strictement
  hyberbolique via l'approximation parabolique.
\newblock {\em C. R. Acad. Sci. Paris S\'er. I Math.}, 319(4):377--382, 1994.

\bibitem{Glimm}
J.~Glimm.
\newblock Solutions in the large for nonlinear hyperbolic systems of equations.
\newblock {\em Comm. Pure Appl. Math.}, 18:697--715, 1965.

\bibitem{GriewankRabier}
A.~Griewank and P.~J. Rabier.
\newblock On the smoothness of convex envelopes.
\newblock {\em Trans. Amer. Math. Soc.}, 322(2):691--709, 1990.

\bibitem{HoldenRisebro}
H.~Holden and N.~H. Risebro.
\newblock {\em Front tracking for hyperbolic conservation laws}, volume 152 of
  {\em Applied Mathematical Sciences}.
\newblock Springer-Verlag, New York, 2002.

\bibitem{JosephLeFloch:ARMA}
K.~T. Joseph and P.~G. LeFloch.
\newblock Boundary layers in weak solutions of hyperbolic conservation laws.
\newblock {\em Arch. Ration. Mech. Anal.}, 147(1):47--88, 1999.

\bibitem{JosephLeFloch}
K.~T. Joseph and P.~G. LeFloch.
\newblock Boundary layers in weak solutions of hyperbolic conservation laws.
  {II}. {S}elf-similar vanishing diffusion limits.
\newblock {\em Commun. Pure Appl. Anal.}, 1(1):51--76, 2002.

\bibitem{JoLeF:nc}
K.~T. Joseph and P.~G. LeFloch.
\newblock {Singular limits for the Riemann problem: general diffusion,
  relaxation, and boundary conditions}.
\newblock In {\em New analytical approach to multidimensional balance laws}. O.
  Rozanova ed., Nova Press, 2007.

\bibitem{kal}
A.~S. Kala{\v{s}}nikov.
\newblock Construction of generalized solutions of quasi-linear equations of
  first order without convexity conditions as limits of solutions of parabolic
  equations with a small parameter.
\newblock {\em Dokl. Akad. Nauk SSSR}, 127:27--30, 1959.

\bibitem{KHass}
A.~Katok and B.~Hasselblatt.
\newblock {\em Introduction to the modern theory of dynamical systems},
  volume~54 of {\em Encyclopedia of Mathematics and its Applications}.
\newblock Cambridge University Press, Cambridge, 1995.
\newblock With a supplementary chapter by Katok and Leonardo Mendoza.

\bibitem{KawShi:normal}
S.~Kawashima and Y.~Shizuta.
\newblock {On the normal form of the symmetric hyperbolic-parabolic systems
  associated with the conservation laws}.
\newblock {\em T\^{o}hoku Math. J.}, 40:449--464, 1988.

\bibitem{lax}
P.~D. Lax.
\newblock Hyperbolic systems of conservation laws. {II}.
\newblock {\em Comm. Pure Appl. Math.}, 10:537--566, 1957.

\bibitem{Liu:TAMS}
T.~P. Liu.
\newblock The {R}iemann problem for general {$2\times 2$} conservation laws.
\newblock {\em Trans. Amer. Math. Soc.}, 199:89--112, 1974.

\bibitem{Liu:rie}
T.~P. Liu.
\newblock The {R}iemann problem for general systems of conservation laws.
\newblock {\em J. Differential Equations}, 18:218--234, 1975.

\bibitem{Liu:adm}
T.~P. Liu.
\newblock The entropy condition and the admissibility of shocks.
\newblock {\em J. Math. Anal. Appl.}, 53(1):78--88, 1976.

\bibitem{MajdaPego}
A.~Majda and R.~L. Pego.
\newblock Stable viscosity matrices for systems of conservation laws.
\newblock {\em J. Differential Equations}, 56(2):229--262, 1985.

\bibitem{Perko}
L.~Perko.
\newblock {\em Differential equations and dynamical systems}, volume~7 of {\em
  Texts in Applied Mathematics}.
\newblock Springer-Verlag, New York, third edition, 2001.

\bibitem{Serre:book}
D.~Serre.
\newblock {\em Systems of conservation laws. 1 and 2}.
\newblock Cambridge University Press, Cambridge, 1999.
\newblock Translated from the 1996 French original by I. N. Sneddon.

\bibitem{Tu1966}
V.~A. Tup{\v{c}}iev.
\newblock The problem of decomposition of an arbitrary discontinuity for a
  system of quasi-linear equations without the convexity condition.
\newblock {\em \u Z. Vy\v cisl. Mat. i Mat. Fiz.}, 6:527--547, 1966.

\bibitem{Tzavaras:JDE}
A.~E. Tzavaras.
\newblock Elastic as limit of viscoelastic response, in a context of
  self-similar viscous limits.
\newblock {\em J. Differential Equations}, 123(1):305--341, 1995.

\bibitem{Tz}
A.~E. Tzavaras.
\newblock Wave interactions and variation estimates for self-similar
  zero-viscosity limits in systems of conservation laws.
\newblock {\em Arch. Rational Mech. Anal.}, 135(1):1--60, 1996.

\end{thebibliography}
\end{document}